\documentclass[onefignum,onetabnum]{siamonline220329}

\usepackage{lipsum}

\usepackage[utf8]{inputenc}
\usepackage{mathrsfs, xcolor,verbatim,bbm,amsmath,amsfonts,amssymb,nicefrac,enumitem,amsxtra}
\usepackage{bm,mathtools,xparse,etoolbox}
\usepackage{bbm}
\usepackage{comment}
\usepackage{graphicx}
\usepackage{epstopdf}
\usepackage[multiple]{footmisc}
\usepackage{algorithm}  
\usepackage{algpseudocode}  

\usepackage{hyperref}

\ifpdf
  \DeclareGraphicsExtensions{.eps,.pdf,.png,.jpg}
\else
  \DeclareGraphicsExtensions{.eps}
\fi

\usepackage{enumitem}
\setlist[enumerate]{leftmargin=.5in}
\setlist[itemize]{leftmargin=.5in}

\newsiamremark{remark}{Remark}
\newsiamremark{hypothesis}{Hypothesis}
\crefname{hypothesis}{Hypothesis}{Hypotheses}
\newsiamthm{claim}{Claim}
\newsiamremark{question}{Question}

\headers{Cluster-based classification with neural ODEs via control}{A. Álvarez, R. Orive, and E. Zuazua}

\title{Cluster-based classification with neural ODEs via control}

\author{Antonio \'Alvarez-L\'opez\thanks{Departamento de Matem\'{a}ticas,
Universidad Aut\'{o}noma de Madrid,
28049 Madrid, Spain \\(\email{antonio.alvarezl@uam.es}).}
\and Rafael Orive-Illera\footnotemark[2] \thanks{Instituto de Ciencias Matem\'aticas, CSIC-UAM-UC3M-UCM, Madrid, Spain
(\email{rafael.orive@icmat.es}).}
\and Enrique Zuazua\footnotemark[2] \thanks{Chair for Dynamics, Control, Machine Learning, and Numerics, Alexander von Humboldt-Professorship, Department of Mathematics,  Friedrich-Alexander-Universit\"at Erlangen-N\"urnberg,
91058 Erlangen, Germany. \\ Chair of Computational Mathematics, Fundaci\'{o}n Deusto, Av. de Universidades 24, 48007 Bilbao, Basque Country, Spain (\email{enrique.zuazua@fau.de}, \url{https://dcn.nat.fau.eu/zuazua/}).}}

\usepackage{amsopn}

\newcommand{\R}{\mathbb{R}}

\newcommand{\N}{\mathbb{N}}

\newcommand{\Z}{\mathbb{Z}}

\newcommand{\cB}{\mathcal{B}}

\newcommand{\cH}{\mathcal{H}}

\newcommand{\cR}{\mathcal{R}}
\newcommand{\cS}{\mathcal{S}}

\newcommand{\cX}{\mathcal{X}}

\newcommand{\bfa}{\mathbf{a}}
\newcommand{\bfb}{\mathbf{b}}
\newcommand{\bfc}{\mathbf{c}}

\newcommand{\bfe}{\mathbf{e}}

\newcommand{\bfg}{\mathbf{g}}

\newcommand{\bfu}{\mathbf{u}}
\newcommand{\bfv}{\mathbf{v}}
\newcommand{\bfw}{\mathbf{w}}
\newcommand{\bfx}{\mathbf{x}}

\newcommand{\scrD}{\mathscr{D}}

\newcommand{\scrX}{\mathscr{X}}
\newcommand{\scrY}{\mathscr{Y}}

\ifpdf
\hypersetup{
  pdftitle={Cluster-based classification with neural ODEs via control},
  pdfauthor={A. Álvarez, R. Orive, and E. Zuazua}
}
\fi

\begin{document}

\maketitle

\begin{abstract}

We address binary classification using neural ordinary differential equations from the perspective of simultaneous control of $N$ data points. We consider a single-neuron architecture with parameters fixed as piecewise constant functions of time. In this setting, the model complexity can be quantified by the number of control switches. Previous work has shown that classification can be achieved using a point-by-point strategy that requires $O(N)$ switches. We propose a new control method that classifies any arbitrary dataset by sequentially steering clusters of $d$ points, thereby reducing the complexity to $O(N/d)$ switches. 
The optimality of this result, particularly in high dimensions, is supported by some numerical experiments. Our complexity bound is sufficient but often conservative because same-class points tend to appear in larger clusters, simplifying classification. This motivates studying the probability distribution of the number of switches required. We introduce a simple control method that imposes a collinearity constraint on the parameters, and analyze a worst-case scenario where both classes have the same size and all points are i.i.d. Our results highlight the benefits of high-dimensional spaces, showing that classification using constant controls becomes more probable as $d$ increases.

\end{abstract}



\section{Introduction}

At the heart of machine learning lies supervised learning \cite{dlnature}, a framework that has been successfully applied in a vast number of domains \cite{imagelecun,languagebengio}. The main objective is to learn an unknown mapping $F:\scrX\to\scrY$. To achieve this, a model $\hat F:\scrX\to\scrY$ to approximate $F$ is constructed by minimizing a loss function, using only the available---possibly noisy---values of $F$ over a finite dataset $\scrD\subset\scrX\times\scrY$. 
Our 
focus is on evaluating the minimal complexity required for $\hat{F}$ to fit the points in $\scrD$ without error.

In the context of data classification, the range of $F$ is finite, and its elements are referred to as labels. Over the years, a wide variety of models have been developed, with notable examples including linear discriminants \cite{fisher}, support vector machines \cite{svm}, random forests \cite{randomforests}, and neural networks \cite{rosenblatt1958perceptron}. 
In \cite{weinan}, a methodology from a control perspective is proposed, based on modeling deep residual networks (ResNets, \cite{he2016residual}) as continuous-time dynamical systems known as neural ordinary differential equations (neural ODEs). 

Neural ODEs have seen the development of several variants \cite{ChenRBD18, haber, massaroli2020dissecting}, yet the standard form remains as
\begin{equation}\label{eq:nodegen}
\left\{
\begin{array}{l}
             \dot \bfx(t)= W(t)\,\boldsymbol{\sigma} (A(t)\,\bfx(t)+\bfb(t)),\hspace{1cm} t\in(0,T),\\
             \bfx(0) = \bfx_0,
\end{array}\right.
\end{equation}
where:
\begin{itemize}
    \item $\bfx_0\in\R^d$ is an input point;
    \item $(W,A,\bfb)\in L^\infty((0,T),\R^{d\times p}\times\R^{p\times d}\times\R^p)$ are parameters to be trained;
    \item $d,p\geq1$ are the state dimension and the width of the model, respectively;
    \item $\boldsymbol{\sigma}:\R^p\to\R^p$ is a prefixed nonlinear Lipschitz function applied component-wise.
\end{itemize} 
Existence and uniqueness of the solution to \cref{eq:nodegen} is guaranteed by the Cauchy–Lipschitz theorem, ensuring the well-definedness of the flow map
\begin{equation}\label{flow}
\Phi_t(\cdot;W,A,\bfb):\bfx_0\in\R^d\longmapsto\bfx(t)\in\R^d,\hspace{1cm}\text{for }t\in[0,T].
\end{equation}
The formulation of \cref{eq:nodegen} naturally frames supervised learning as a control problem. Here, the input space is $\mathscr{X}=\R^d$, and the parameters $(W,A,\bfb)$ serve as controls that simultaneously guide all input points toward their respective target positions in $\R^d$. To match the output space, a mapping $g:\R^d\to\mathscr{Y}$ is introduced as a final layer. The complete model is thus defined by the composition $\hat F = g\circ \Phi_T$.

We focus on binary classification with a hard classifier, whereby $\mathscr{Y}=\{1,0\}$ and $g$ is the characteristic function of a fixed set. Nonetheless, our results can be extended to any multiclass setting by fixing $g$ as a weighted sum of predefined characteristic functions, each corresponding to a distinct label.

Neural ODEs were originally conceived as a tool for understanding deep ResNets, but they have since made a significant impact on machine learning. Their continuous-time framework facilitates mathematical analysis 
and provides practical benefits like incorporating structure or the design of new discrete schemes. For more details, we refer to \cite{weinan}.

\subsection*{Notation}
\begin{itemize}
\item Scalars are denoted by plain letters, vectors by boldface letters, and matrices by uppercase letters. The scalar product of two vectors $\bfu$, $\bfv$ is written as $\bfu\cdot\bfv$.
\item Subscripts identify elements of a set. Superscripts identify components of a vector.
   \item $\{\bfe_1,\dots,\bfe_d\}$ denotes the canonical basis in $\R^d$.
   \item $\mathbb{S}^{d-1}$ denotes the $(d-1)-$dimensional sphere in $\R^d$.
    \item The cardinality of a set $\cX$ is denoted by $|\cX|$.
   \item For $x\in\R$, we write $\lceil x\rceil\coloneqq\min\left\{n\in\Z:n\geq x\right\}$, $\lfloor x\rfloor\coloneqq\max\left\{n\in\Z:n\leq x\right\}.$
\end{itemize}

\section{Problem formulation and main results}\label{sec:2}
Let $\scrD\coloneqq\{(\bfx_n,y_n)\}\subset\R^d\times\{1,0\}$ be a finite dataset such that $\bfx_i\not=\bfx_j$ for all $i\neq j$. We define the classes $\cR$ (red circles) and $\cB$ (blue crosses) by
\begin{equation}\label{eq:RB}
\cR=\{\bfx_n\in\R^d:\; (\bfx_n,1)\in\scrD\},\hspace{1cm}
 \cB=\{\bfx_n\in\R^d:\; (\bfx_n,0)\in\scrD\}.
\end{equation}
We adopt the simplified version of neural ODEs with one-neuron width.  Namely, we set $p=1$ in \cref{eq:nodegen}, which yields
 \begin{equation}\label{eq:node}
    \dot \bfx(t) = \bfw(t)\sigma(\bfa(t)\cdot\bfx(t)+b(t)),\hspace{1cm} \text{for }t\in(0,T).
\end{equation} 
Here, $\sigma(\cdot)=(\cdot)_+$ is the rectified linear unit (ReLU), while $\theta=(\bfw,\bfa,b)$ belongs to \begin{equation*}
    \Theta_T\coloneqq L^\infty\left((0,T)\;;\;\mathbb{S}^{d-1}\times\R^d\times \R\right).
 \end{equation*}
and is assumed to be piecewise constant. Imposing this constraint reduces optimization to a finite-dimensional space \cite{massaroli2020dissecting}, while also inducing a layered structure akin to that of a discrete ResNet  \cite{domzuazua,ALVAREZLOPEZ2024106640}. For any control $\theta \in \Theta_T$, we define the complexity of  \cref{eq:node} as the number of finite-jump discontinuities (or switches) of $\theta$ over $(0,T)$, denoted by $L$.

Classification can essentially be interpreted as transforming data into representations in which different classes are separable. In the neural ODE framework, given a dataset $(\cR,\cB)$, the goal is to find a control $\theta\in\Theta_T$ for  \cref{eq:node} that induces a finite-time flow map $\Phi_T$  satisfying
\begin{equation*}
     \Phi_T(\cR;\theta)\subset \tau_\cR\hspace{1cm}\text{and}\hspace{1cm}\Phi_T(\cB;\theta)\subset \tau_\cB,
 \end{equation*}
 where $(\tau_\cR,\tau_\cB)$ is a pair of linearly separable regions of $\R^d$
. 
For simplicity, we fix $(\tau_\cR,\tau_\cB)$ to be half-spaces of the form
\begin{equation}\label{eq:targets}
 \left(\{x^{(i)}>1\},\{ x^{(i)}\leq1\}\right)\quad\text{or}\quad \left(\{x^{(i)}\leq1\},\{ x^{(i)}>1\}\right)\hspace{1cm}\text{for }i\in\{1,\dots,d\},
 \end{equation}
and the hyperplane $\{x^{(i)}=1\}$ as decision boundary. We make use of the dynamics of \cref{eq:node} in a constructive manner, by carefully defining each value of the piecewise constant control $\theta\in\Theta_T$. On the layer $t\in(t_{k-1},t_k)\subset(0,T)$, the parameters $\bfa(t)\equiv \bfa\in\R^{d}$ and $b(t)\equiv b\in\R$ determine the hyperplane $H:\bfa\cdot\bfx+b=0$ and the two half-spaces
\begin{equation*}
H^+\coloneqq\big\{\bfx\in\R^d:\bfa\cdot\bfx+b>0\big\},\hspace{1cm} H^-\coloneqq\R^d\setminus H^+.
\end{equation*}
The half-space $H^-$ remains fixed because $\sigma(\bfa\cdot\bfx + b)=0$ for all $\bfx\in H^-$. in contrast, each point $\bfx\in H^+$ evolves according to the vector field $\bfw(\bfa\cdot\bfx+b)$. The direction of this field is constant and determined by $\bfw\in\mathbb{S}^{d-1}$, while its magnitude at each $\bfx\in H^+$ is equal to the distance from $\bfx$ to the hyperplane $H$.

\begin{figure}[t!]  
\centering
    \begin{minipage}[b]{0.55\textwidth}
        \centering
        \includegraphics[scale=0.07]{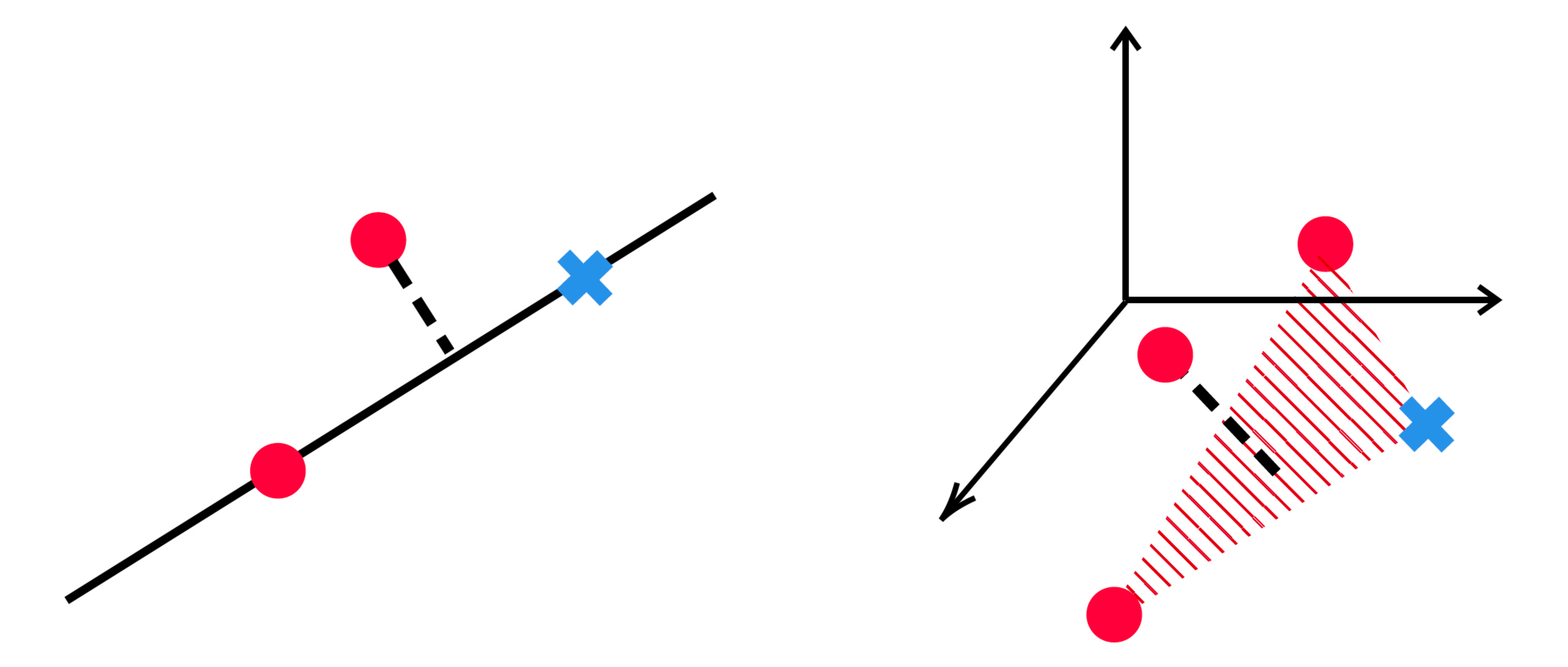}
        \caption{$\cX\subset\R^2$ is in general position if no three points of $\cX$ lie on the same line. $\cX\subset\R^3$ is in general position if, additionally, no four points lie on the same plane.}
        \label{fig:genpos}
    \end{minipage}
      \hfill
       \begin{minipage}[b]{0.42\textwidth}
        \centering
\includegraphics[scale=0.35]{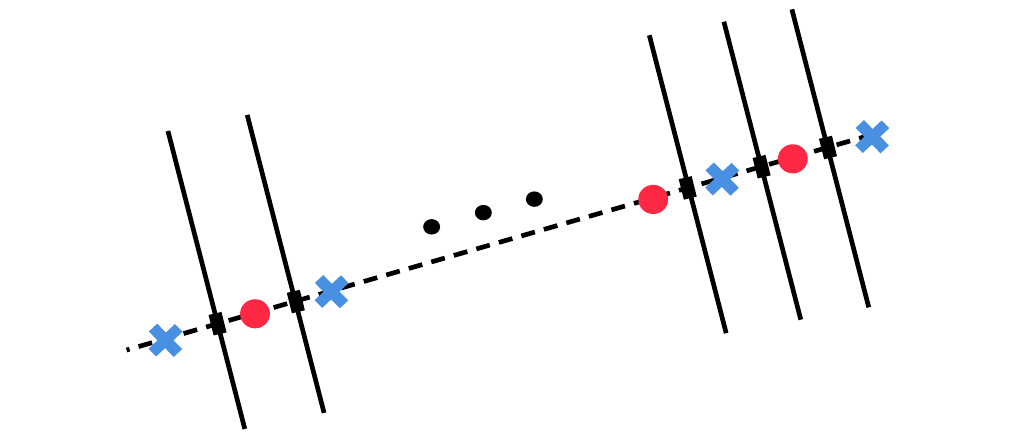}
\caption{Configuration where the maximum number of $2N-1$ hyperplanes is required to separate the points.}
\label{fig:Z2N-1} 
    \end{minipage}
\end{figure}

This overview of the dynamics suggests that classifying data points using \cref{eq:node} is intrinsically related to identifying a set of hyperplanes that separate them by labels. All in all, there is a clear relationship between the complexity of the controls in \cref{eq:node} and the geometric structure of the data distribution. The following questions naturally arise:

    \begin{question}\label{quest1}
   What is the minimum number of discontinuities required to ensure that any dataset with a fixed number of points can be classified?
\end{question}
\begin{question}\label{quest2}
   What is the probability that a given dataset can be classified using exactly $L=k$ discontinuities, for any $k\geq0$?
    \end{question}
Regarding \cref{quest1}, it was established in \cite{domzuazua} that $L=3\min\left\{|\cR|,|\cB|\right\}$ switches are sufficient to classify $(\cR,\cB)$ with \cref{eq:node}. The proof involves sequentially steering each point in $\cR\cup \cB$ individually to its corresponding target $\tau_\cR$ or $\tau_\cB$. Our first main result refines this bound by relying on the concept of general position, a classical notion in combinatorial geometry \cite{mumford1965geometric, cover1965,cover_number_1967,sontag_shattering_1997}. We consider the following definition, illustrated in \cref{fig:genpos}:

    \begin{definition}\label{def:gp}
  A set $\cX\subset\R^d$ is in \emph{general position} if, for every $0\leq k \leq d-1$, no affine subspace in $\R^d$ of dimension $k$ contains more than $k+1$ points of $\cX$. 
\end{definition}
Assuming this mild condition---easily achieved by slight perturbations of the dataset---we construct a family of pairwise parallel hyperplanes that enclose all points of one class by subsets of size $d$. We can then define a simultaneous control method that classifies the points within these subsets rather than individually, thereby reducing the bound for $L$:

\begin{theorem}\label{thm1}
    Let $d\geq 2$. For any dataset $(\cR,\cB)$ defined as in \cref{eq:RB} in general position and any pair of target sets $(\tau_\cR,\tau_\cB)$ defined as in \cref{eq:targets}, there exist $T>0$ and a piecewise constant control
$\theta\in\Theta_T$ whose number of discontinuities is  \begin{equation*}
    L=4\left\lceil \frac{\min\left\{|\cR|,|\cB|\right\}}{d}\right\rceil-1,
\end{equation*}
such that the flow map of the neural ODE \cref{eq:node} satisfies    $\Phi_T(\cR;\theta)\subset \tau_\cR$ and $\Phi_T(\cB;\theta)\subset \tau_\cB$.
\end{theorem}
\cref{thm1} holds almost surely as long as the points are sampled from a non-singular measure. However, numerous data points sharing the same label may initially be clustered together, enabling their control with fewer parameters. This observation motivates the introduction of a probabilistic framework to address \cref{quest2}.  

We consider the worst-case scenario where all points are i.i.d., and the sizes are fixed at $|\cR|=|\cB|=N$. 
To simplify the analysis, we introduce a new control method that restricts $\bfa(t)$ to be a constant, although optimally chosen. This constraint will allow us to determine the exact probability distribution of $L$ for this method.


\begin{theorem}\label{thm2} Let $d\geq 2$ and $N\geq1$. Consider any dataset $(\cR,\cB)$ defined as in  \cref{eq:RB}, with $|\cR|=|\cB|=N$. Assume every point $\bfx\in\cR\cup\cB$ is independently sampled from an absolutely continuous probability measure on $[0,1]^d$. 
Then, with probability 1, there exist $T>0$, a pair of target sets $(\tau_\cR,\tau_\cB)$ defined as in \cref{eq:targets}, and $\theta=(\bfw,\bfa,b)\in\Theta_T$ with \begin{equation*}\bfa(t)\in\{\bfe_1,\dots,\bfe_d\}\hspace{1cm}\text{constant for all  }t\in(0,T),\end{equation*}
such that the flow map of the neural ODE \cref{eq:node} satisfies 
     $\Phi_T(\cR;\theta)\subset \tau_\cR$ and $\Phi_T(\cB;\theta)\subset \tau_\cB$.
     
Furthermore, $\theta$ is piecewise constant with $L$ discontinuities following, for $0\leq k\leq2N-2$,
\begin{equation}\label{eq:distribL}
    \mathbb{P}(L\geq k)=\left(\sum_{p=\lceil\frac{k}{2}+1\rceil}^N\binom{N-1}{p-1}^2+\sum_{p=\lceil\frac{k+1}{2}\rceil}^{N-1}\binom{N-1}{p}\binom{N-1}{p-1}\right)^d2^d\binom{2N}{N}^{-d}.
\end{equation}

\end{theorem}
\begin{remark}
    Stirling's approximation simplifies \cref{eq:distribL} for large values of $d$ and $N$. If we set $k=1$, the formula reduces to
\begin{equation*}
\mathbb{P}(L\geq 1) = \left(1-2\binom{2N}{N}^{-1}\right)^d \sim \exp\left\{\frac{-\sqrt{\pi N}}{2^{2N-1}}\;d\right\}.   
\end{equation*}
Assuming both $d$ and $N$ grow according to some relation $d = d(N)$, we can deduce that
\begin{equation*}  \text{if}\hspace{1cm}\lim_{N\rightarrow\infty}\frac{2^{2N}}{d(N)\sqrt{N}}=0 \hspace{1cm}\text{then}\hspace{1cm}\lim_{N\rightarrow\infty}\mathbb{P}(L=0)=1. \end{equation*}
This observation reveals an explicit trade-off between $d$ and $N$ that allows classifying all points employing  a constant control, or equivalently, with an autonomous neural ODE.
\end{remark}

The control method of \cref{thm2} is designed so that distribution \cref{eq:distribL} is derived from the number of hyperplanes required to separate the points by labels. The maximum value $L=2N-2$, corresponding to $2N-1$ hyperplanes, also improves upon the bound of \cite{domzuazua}.

In the following theorem, we take a geometric perspective on the problem to show that the maximum number of hyperplanes required to separate the two classes is indeed $2N-1$. Moreover, we characterize the pathological point configurations that attain this maximum, as illustrated in \cref{fig:Z2N-1}.

\begin{theorem}\label{thm3}
Let $d,N\geq 1$. For any dataset $(\cR,\cB)$ defined as in \cref{eq:RB} with $|\cR|=|\cB|=N$, the maximum number of hyperplanes in $\R^d$ required to separate the points by labels is $2N-~1$. Furthermore, this maximum is attained if and only if the points of $\cR$ and $\cB$ are collinear and alternating, i.e.,  there exist $\bfu,\bfu_0\in\R^d$ and $-\infty<\lambda_1<\cdots<\lambda_{2N}<+\infty$  such that
\begin{equation*}
\bfu_0+\lambda_{2k-1}\bfu\in\cR\hspace{1cm}\text{and}\hspace{1cm}\bfu_0+\lambda_{2k}\bfu\in\cB\hspace{1cm}\text{for all }k\in\{1,\dots,N\}.
\end{equation*}
\end{theorem}
We note that deriving the maximum is relatively straightforward; the main challenge lies in showing that it is attained only when all points are collinear and alternating.


\subsection*{Roadmap}

In \cref{sec:3}, we develop the mathematical framework, which mainly relies on hyperplane separability. We begin by proving \cref{thm3} and then present combinatorial results---first for $d=1$ and later in higher dimensions---that culminate in \cref{thm2}. In \cref{sec:4}, assuming all points are in general position (see \cref{def:gp}), we construct a family of pairwise parallel hyperplanes that separate the dataset, with each pair enclosing exactly $d$ points of the same class. Next, we prove \cref{thm1} by defining a classification method that sequentially controls these subsets of size $d$, and analyze alternative activation functions that further reduce the value of $L$. Although most results are static and not specific to neural ODEs, \cref{thm1,thm2} incorporate dynamic algorithms, distinguishing them from standard linear classifiers. In \cref{sec:5}, we perform a computational test using gradient-based training to estimate the minimal complexity required for neural ODEs to classify datasets of fixed size. We then compare these results to the complexity obtained in \cref{thm1}. We conclude with a summary of our main contributions, a discussion on connections and extensions, and open questions for future work.

\subsection*{Related work}
The theoretical study of neural ODEs from a control theory perspective has gained significant attention in recent years, opening several and promising research directions.

One line of research explores the approximation power of neural ODE flows, using either geometric techniques based on Lie brackets \cite{TaGh, CuLaTe, EZOP} or constructive methods for simultaneous control \cite{Qianxiao2022, domzuazua}. In \cite{cheng2023interpolation}, the authors conducted a detailed study of controllability and its connection with density in $L^p$ spaces for $p<+\infty$. Universal approximation in $L^\infty$ was established for a specific class of diffeomorphisms in \cite{ishikawa_universal_2023}. We also highlight the work of \cite{uat1neuronresnet}, which demonstrates density in $L^1(\R^d)$ for the family of ResNets with one neuron per hidden layer, corresponding to the forward Euler discretization of \cref{eq:node}.

Another approach studies the training process as an optimal control problem. In \cite{e_mean-field_2018}, population risk minimization in deep learning is formulated as a mean-field optimal control problem, and Hamilton–Jacobi–Bellman and Pontryagin optimality conditions are derived. \cite{qianxiaomaxprinciple} proposes a modification of the successive approximations method by augmenting the Hamiltonian to solve Pontryagin’s maximum principle, resulting in an alternative training algorithm with rigorous error estimates. \cite{bonnet_measure_2023} employs a continuity equation to study the mean-field dynamics, examining the existence and uniqueness of minimizers in the optimal control problem with $L^2$-regularization and establishing a mean-field maximum principle.  In \cite{24isobejml}, the result is improved by introducing a kinetic regularization term in the loss function, and proving the existence of minimizers. \cite{Geshkovski_Zuazua_2022} considers the classical empirical risk minimization, deriving a manifestation of the turnpike property through specific regularization terms. Additionally, \cite{EGe} reduces the complexity for data classification by introducing an $L^1$-norm penalty, which promotes temporal sparsity in the control. 

We follow a third direction that consists of estimating the complexity required to control  $N$ points.  To the best of our knowledge, the only study addressing this problem in the one-neuron width model \cref{eq:node} is \cite{domzuazua}, where the bound of $O(N)$ switches is established. We build upon their results by reducing to $O(N/d)$ switches, and deriving a probabilistic result.  In \cite{ALVAREZLOPEZ2024106640}, the result of \cite{domzuazua} is generalized to any finite width $p$, resulting in a complexity of $O(N/p)$. The problem is also explored in \cite{momentum}, which examines a neural ODE with a second-order time derivative to further enrich the dynamics, maintaining a complexity of $O(N)$. Furthermore, \cite{domzuazuaNormFlows} investigates the controllability of probability measures for the continuity equation that extends \cref{eq:node} and considering  errors in total variation.

\section{Classification via canonical separability}\label{sec:3}

Suppose that all points of the dataset \cref{eq:RB} are independently sampled from a particular absolutely continuous probability measure on $[0,1]^d$, with $|\cR|=|\cB|=N$ fixed, 
and that they satisfy the following condition:
\begin{equation}\label{eq:condcansep}     x_n^{(j)}\neq x_m^{(j)}, \quad \hbox{ for all } j\in\{1,\dots,d\}\quad\hbox{and}\quad n\neq m. \end{equation}
Note that \cref{eq:condcansep} is fulfilled with probability 1. We now introduce our concept of separability:

\begin{definition}
A finite family $\cH$ of affine $(d-1)$-dimensional hyperplanes in $\R^d$ \emph{separates} $(\cR,\cB)$ if they break the cube $[0,1]^d$ into polyhedra, each of them containing points of at most one of the two sets. We say that such $\cH$ is a collection of \emph{separating hyperplanes} for $(\cR,\cB)$.  
\end{definition}

Note that some polyhedra might not contain any point from the two sets. Now, for any pair $(\cR,\cB)$ as in \cref{eq:RB}, and any collection of hyperplanes $\cH$ in $\R^d$, let us consider
 \begin{equation}\label{eq:ZdN}
   Z_{d,N}(\cR,\cB)\coloneqq\min\left\{|\cH|\,:\,\cH\text{ separates }(\cR,\cB)\right\}.  
 \end{equation}
As defined, $Z_{d,N}$ is a random variable that maps each pair $(\cR,\cB)$ to the minimum cardinality of any collection of separating hyperplanes for $(\cR,\cB)$.

Theorem~\ref{thm2} is based on separating  $(\cR,\cB)$ using hyperplanes that are orthogonal to an optimally chosen canonical vector.  We refer to this approach as \emph{canonical $k$-separability}, drawing on the concept of $k$-\emph{separability} introduced in \cite{duchksep}. A finite set $\{\bfx_n\}\subset\R^d$ with binary labels is $k$-separable if there exists $\bfw \in \mathbb{S}^{d-1}$ such that the projections $\bfw \cdot \bfx_n$ can be divided into $k$ disjoint intervals, each containing only projections with the same label. 

Canonical $k$-separability is a similar concept, but it constrains $\bfw$ to the canonical basis. We will determine the exact probability that $(\cR,\cB)$ is $k$-separable under this constraint, which requires adapting the definition of $Z_{d,N}$ in \cref{eq:ZdN}.

For each $i=1,\dots,d$, the usual projection $\pi^i:\R^d\to \R$ is defined by $\pi^{i}(\bfx)=x^{(i)}$ for all $\bfx\in\R^d$. Let $\Pi^{i}$ be the pointwise extension of $\pi^i$, given by
\begin{equation*}
\Pi^{i}(\cX)=\left\{\pi^i(\bfx_1),\dots,\pi^i(\bfx_N)\right\}=\left\{x_1^{(i)},\dots,x_N^{(i)}\right\}
\end{equation*}
for all $\cX=\{\bfx_n\}_{n=1}^N\subset\R^d$ satisfying condition \cref{eq:condcansep}, which ensures $|\Pi^{i}(\cX)|=N$. 

\begin{figure}[t!]  
\centering
\includegraphics[width=0.48\textwidth]{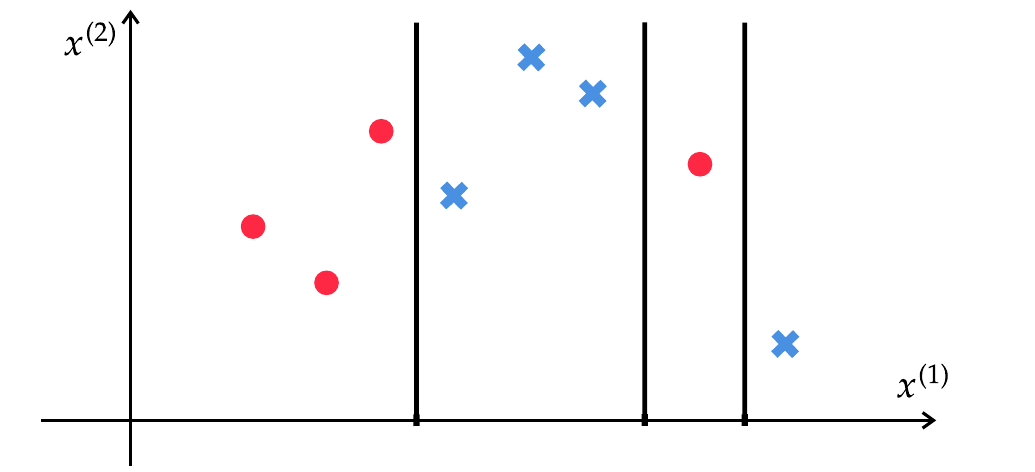}
\includegraphics[width=0.48\textwidth]{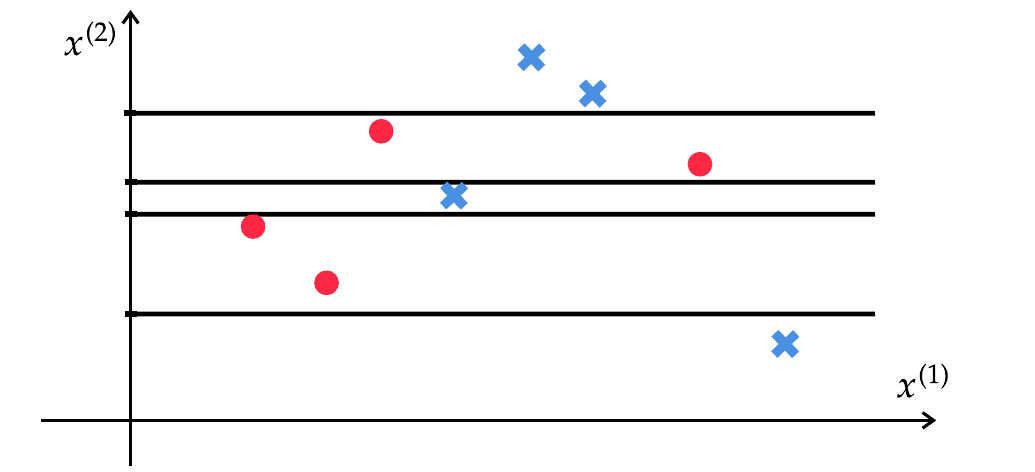}
\caption{$Z_{d,N}^1=3$ and $Z_{d,N}^2=4$ computed by projecting the data on the respective axes $x^{(1)}$ and $x^{(2)}$.}
\label{fig:projections}
\end{figure}

We introduce the random variable $Z_{d,N}^i$ such that
\begin{equation}\label{eq:ZdNi}
Z_{d,N}^i(\cR,\cB)=Z_{1,N}\left(\Pi^i(\cR),\Pi^i(\cB)\right)
\end{equation}
for any pair $(\cR,\cB)$ defined as in \cref{eq:RB}. The value $Z_{d,N}^i(\cR,\cB)$ represents the minimum number of points required to separate the projections of both sets on the $i$-th Cartesian axis. Those points determine a family of parallel hyperplanes in $\R^d$ given by the equations $x^{(i)}=H_j$  with $j=1,\dots, Z_{d,N}^i(\cR,\cB)$, that separate $(\cR,\cB)$ in $\R^d$, see \cref{fig:projections}.

Given $d,N\geq1$, we define 
\begin{equation}\label{eq:ZdNorth}
Z^\perp_{d,N}\coloneqq\min\left\{Z_{d,N}^1,\dots,Z_{d,N}^d\right\}.
\end{equation}
The variable $Z_{d,N}^\perp$ represents the minimum cardinality of any family of separating hyperplanes, all of which are perpendicular to some Cartesian axis.
For the example shown in \cref{fig:projections}, we would compute $Z_{d,N}^\perp(\cR,\cB)=3$.

\smallskip

Now, we prove \cref{thm3} to determine and characterize the maximum value of $Z_{d,N}$:

\begin{proof}[Proof of \cref{thm3}]  
By relabeling the points, with no loss of generality we can assume 
\begin{equation*}
    \cR = \{\bfx_n\}_{1\leq n\leq N}\hspace{1cm}\text{and}\hspace{1cm} \cB = \{\bfx_{N+n}\}_{1\leq n\leq N}.
\end{equation*}
First, we prove
\begin{equation}\label{eq:maxnumd}
\max_{(\cR,\cB)}Z_{d,N}(\cR,\cB)=2N-1.
\end{equation}
\noindent {\bf Case $d=1$}. 
When the points of the pair $(\cR,\cB)$ are alternating, i.e.,
\begin{equation}\label{eq:intersp1d}
 x_1<x_{N+1}<x_2<\cdots<x_N<x_{2N}\quad\text{or}\quad x_{N+1}<x_1<x_{N+2}<\cdots<x_{2N}<x_N,
\end{equation}
we have  $Z_{1,N}(\cR,\cB)=2N-1$. Conversely, if the points of $\cR\cup \cB$ are not alternating, then the number of line segments connecting consecutive points with different labels is less than $2N-1$ and prove  \cref{eq:maxnumd} with $d=1$.

\smallskip
\noindent{\bf Case $d>1$}. Since $Z_{d,N}\leq Z^\perp_{d,N}$ as defined in \cref{eq:ZdNorth}, we deduce that $\max Z_{d,N}\leq 2N-1$  using \cref{eq:maxnumd} for $d=1$. To show that $\max Z_{d,N}\geq 2N-1$, assume the $2N$ points of $\cR\cup \cB$ are collinear and alternating. To separate $(\cR,\cB)$, each of the $2N-1$ line segments connecting two consecutive points of different labels must then be intersected by a hyperplane. Since any hyperplane that does not contain the entire line can intersect it in at most one point, we need at least $2N-1$ hyperplanes to separate $(\cR,\cB)$, as illustrated in \cref{fig:Z2N-1}.



\smallskip
\noindent\textbf{Necessity of being collinear and alternating (only if).} If the points of $\cR\cup \cB$ are collinear but not alternating, then the number of line segments connecting consecutive points with different labels is less than $2N-1$. Consequently, fewer than $2N-1$ hyperplanes suffice to separate them, that is, $Z_{d,N}<2N-1$.

To show that collinearity is also a necessary condition to attain the maximum, we proceed by induction on $N$. For $N=1$ the claim is trivial. For $N=2$, if the four points are not collinear, their convex hull is at most a triangle, quadrilateral, or tetrahedron, which implies $Z_{d,N}\leq 2$.

Let $N>2$ and suppose that $\cR\cup\cB$ is not contained in any line, yet $Z_{d,N}(\cR,\cB)=2N-1$. Since the dataset is finite, we can choose a line $r\subset\R^d$ with direction $\bfv\in\mathbb{S}^{d-1}$ such that the orthogonal projection of $\cR\cup\cB$ onto $r$ is injective. The number of points of $r$ required to separate the projected sets is then $2N-1$. Indeed, it is at most $2N-1$ (as per \eqref{eq:maxnumd} for $d=1$), and on the other hand, these points correspond to a family of hyperplanes orthogonal to $r$ that separate $(\cR,\cB)$, so it must be at least $Z_{d,N}(\cR,\cB)=2N-1$. All in all, we can choose hyperplanes $H_1,\dots,H_{2N-1}$ orthogonal to $r$ that separate $(\cR,\cB)$, and moreover, the projected points $\{\bfv\cdot\bfx_i\}_{i=1}^{2N}\subset\R$ must be alternating as in \cref{eq:intersp1d}. Without loss of generality, we may assume
\begin{equation*}
\bfv\cdot \bfx_1 < h_1< \bfv\cdot \bfx_{N+1}< \cdots< \bfv\cdot \bfx_N< h_{2N-1}< \bfv\cdot \bfx_{2N},\end{equation*}
where $h_i = \bfv\cdot \mathbf{h}_i$ and $\mathbf{h}_i = H_i \cap r \subset \R^d$, as illustrated in \cref{fig:1demopropmax} (left). We then define $\cR' = \cR\setminus\{\bfx_1\}$ and $\cB' = \cB\setminus\{\bfx_{N+1}\}$, each containing $N-1$ points.

Now, we show that $Z_{d,N-1}(\cR',\cB')=2N-3$. Suppose that $k$ hyperplanes $G_1,\dots,G_k$ in $\R^d$ suffice to separate $(\cR',\cB')$ for some $1\leq k\leq 2N-4$. Then the pair $(\cR,\cB)$ can be separated using $k+2$ hyperplanes $H_1,H_2,G_1,\dots, G_k$, as illustrated in \cref{fig:1demopropmax} (right). Since $k+2\leq2N-2<2N-1$, this contradicts the fact that $Z_{d,N}(\cR,\cB)=2N-1$.

\begin{figure}[t]  
\centering
\includegraphics[width=0.48\textwidth]{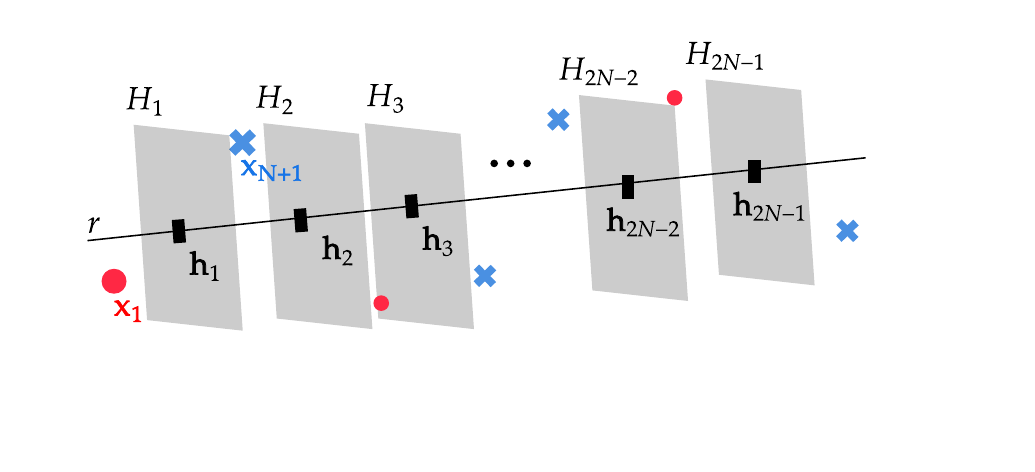}
\includegraphics[width=0.48\textwidth]{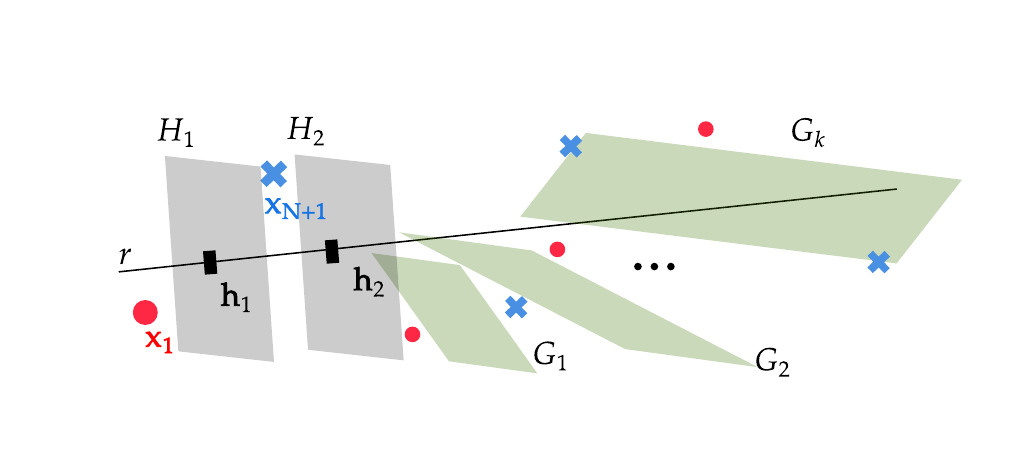}
\caption{Figures supporting the argument presented in the proof of \cref{thm3}.}
\label{fig:1demopropmax}
\end{figure}

\begin{figure}[t]  
\centering
\includegraphics[width=0.48\textwidth]{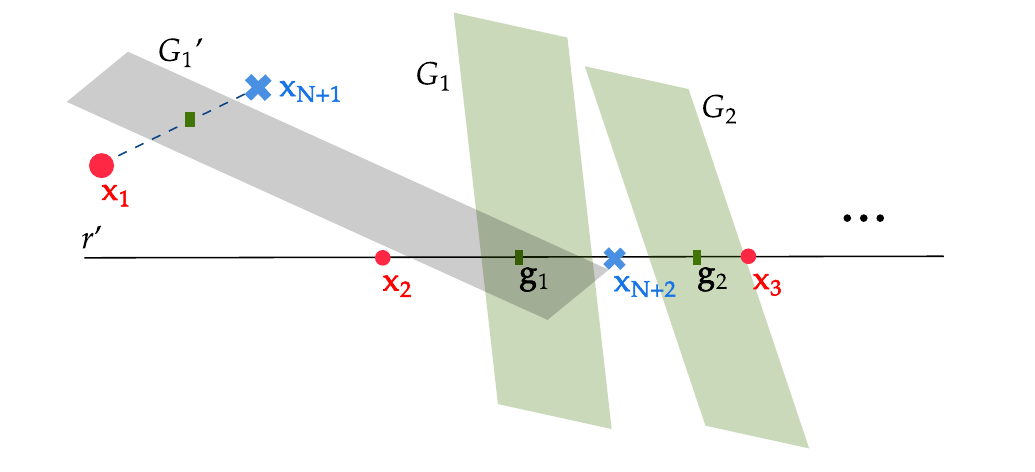}
\includegraphics[width=0.48\textwidth]{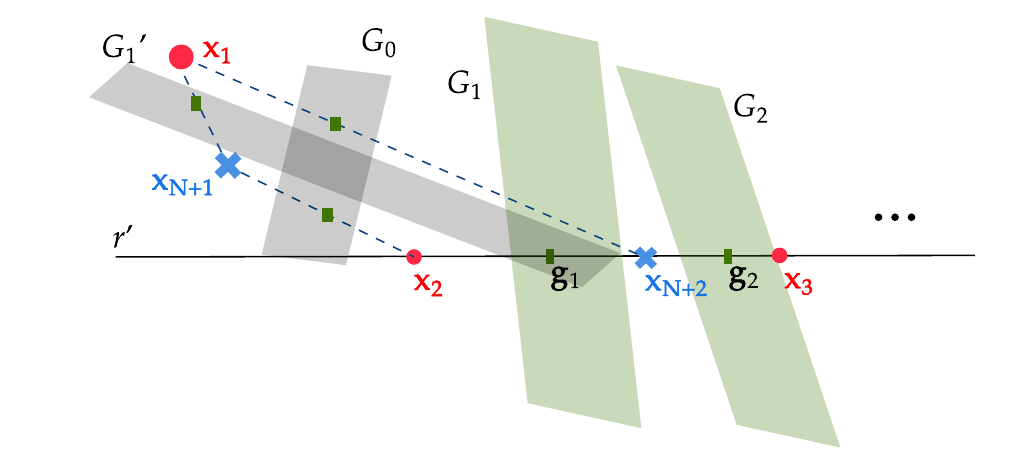}
\caption{Figures supporting the argument presented in the proof of \cref{thm3}.}
\label{fig:2demopropmax}
\end{figure}

The remainder of the proof consists of verifying that $(\cR,\cB)$  can be separated by at most $2N-2$ hyperplanes, which would lead to a contradiction. Since $Z_{d,N-1}(\cR',\cB')=2N-3$, the inductive hypothesis on $N-1$ ensures that $\cR'\cup\cB'$ must be collinear along some line $r' \subset \mathbb{R}^d$.  By assumption, the points in $\cR\cup \cB$ are not collinear, which implies that  $\bfx_1\not\in r'$ or $\bfx_{N+1}\not\in r'$. Let $G_1,\dots,G_{2N-3}$ be hyperplanes that separate $(\cR',\cB')$ and each intersecting $r'$ transversely, that is, we can set $\bfg_j= G_j\cap r'\in\R^d$. We distinguish two cases:
\smallskip\newline
1. If $\bfx_1$, $\bfx_{N+1}$ and $\bfg_1$ are not collinear, we can consider any hyperplane $G_1'$ that intersects the open segment $[\bfx_1,\bfx_{N+1}]=\{t\bfx_1+(1-t)\bfx_{N+1}\,|\,0<t<1\}$ transversely, and which also intersects the line $r'$ transversely at $\bfg_1$. If $\bfx_1$ is on the same side of $G_1'$ as $\bfx_2$, the family $\{G_1',G_2,\dots,G_{2N-3}\}$ separates $(\cR,\cB)$---as illustrated in \cref{fig:2demopropmax} (left)---and $Z_{d,N}(\cR,\cB)=2N-3$ (contradiction). Otherwise, we also consider any hyperplane $G_0$ that intersects transversely with the open segment $[\bfx_2,\bfx_{N+1}]$ and with the open segment $[\bfx_1,\bfx_{N+2}]$. Then, $\{G_0,G_1',G_2,\dots,G_{2N-3}\}$ separates $(\cR,\cB)$, as represented in \cref{fig:2demopropmax} (right), and $Z_{d,N}(\cR,\cB)=2N-2$ (contradiction).\smallskip\newline
2. If $\bfx_1$, $\bfx_{N+1}$ and $\bfg_1$ are collinear,  we can slightly translate $G_1$ and perturb its intersection point with $r'$ from $\bfg_1$ to another point $\bfg_1'$ of the open segment $[\bfx_2,\bfx_{N+2}]$. Thus, we can ensure that $\bfx_1$ and $\bfx_{N+1}$ and $\bfg_1'$ are not collinear and apply case 1.

\end{proof}

To obtain the probability distribution of $Z_{d,N}^\perp$, we first solve the case $d=1$ employing combinatorial techniques:

\begin{lemma}\label{prop:3}
For $N\geq1$, let $Z_{1,N}$ be as defined in \cref{eq:ZdN}. For any $1\leq k\leq 2N-1$, we have
\begin{equation}\label{eq:distrib1d}
\mathbb{P}(Z_{1,N}=k)= \begin{cases}\displaystyle
2\binom{N-1}{p-1}^2 \binom{2N}{N}^{-1} \quad &\text{if} \quad k=2p-1, \vspace{3mm}\\ 
\displaystyle
 2\binom{N-1}{p}\binom{N-1}{p-1} \binom{2N}{N}^{-1}  \quad&\text{if} \quad k=2p.
\end{cases}
\end{equation}
\end{lemma}

\begin{proof}
Since all points in $\cR \cup \cB$ are i.i.d., each possible ordering is equally likely. Therefore, to determine $\mathbb{P}(Z_{1,N} = k)$,  we count the fraction of orderings that have exactly $k$ gaps between consecutive subsets of points with different labels. We compute the number of such favorable configurations by dividing them into two cases based on the parity of $k$: 
\smallskip

\textbf{Case 1:} $k=2p-1$, for $1\leq p\leq N$.
   Any ordering of the points in  $\cR\cup \cB$ that leads to $Z_{1,N}(\cR,\cB)=2p-1$ can be represented as the union of two partitions of $\cR$ and $\cB$, respectively constituted of subsets $R_i\subset \cR$ and $B_i\subset \cB$, for $i\in1,\dots,p$, with lengths $|R_i|=r_i$ and $|B_i|=b_i$,
that satisfy 
\begin{equation}\label{eq:part}
\sum_{i=1}^pr_i=\sum_{i=1}^pb_i=N.\end{equation}
There are $\binom{N-1}{p-1}$ possible partitions $\{R_1,\dots,R_p\}$ of $\cR$,  determined by the choice of $p-1$ gaps between consecutive elements of $\cR$, among the total $N-1$ possibilities. Analogously, there are $\binom{N-1}{p-1}$ possible partitions $\{B_1,\dots,B_p\}$ of $\cB$, to insert in the $p$ gaps and therefore obtain $Z_{1,N}(\cR,\cB)=k$. So there are $\binom{N-1}{p-1}^2$ possible configurations. On the other hand, we can repeat the argument but now considering the $\binom{N-1}{p-1}$ possible partitions of $\cB$ and then inserting $R_i$ into the $p$ gaps. Both situations are symmetric; the only difference between them lies in whether the smallest point of $\cR\cup \cB$ (and thus the first subset of the two partitions) belongs to $\cR$ or $\cB$, as shown in \cref{fig:k12}. 
Consequently, the total number of configurations for $(\cR,\cB)$ that yield $Z_{1,N}(\cR,\cB)=2p-1$ is $2\binom{N-1}{p-1}^2$.\newline
\textbf{Case 2:} $k=2p$, for $1\leq p\leq N-1$. To obtain $Z_{1,N}(\cR,\cB)=2p$, we need an even number of gaps between 
subsets of the same color, so now we must consider the partitions of $\cR$ and $\cB$ that are respectively constituted of subsets $R_i,R_{p+1}\subset \cR$ and $B_i\subset \cB$ for $i\in1,\dots,p$ (or the reverse situation) with lengths $|R_i|=r_i$ 
and $|B_i|=b_i$,
that satisfy 
\begin{equation}\label{eq:part2}
\sum_{i=1}^{p+1}r_i=\sum_{i=1}^pb_i=N.\end{equation}
Like in case 1, there exist $\binom{N-1}{p}$ possible partitions $\{R_1,\dots,R_{p+1}\}$ of $\cR$ and $\binom{N-1}{p-1}$ possible partitions $\{B_1,\dots, B_p\}$ of $\cB$. So there are $\binom{N-1}{p}\binom{N-1}{p-1}$ possibilities, each one depending on a choice of $(r_1,\dots,r_{p+1})\in\N^{p+1}$ and $(b_1,\dots,b_p)\in\N^p$ satisfying \cref{eq:part2}. We take also into account the symmetric situation, when the partitions of $\cR$ and $\cB$ are respectively of the form $\{R_1,\dots,R_p\}$ and $\{B_1,\dots,B_{p+1}\}$, 
 shown in \cref{fig:k34}. Therefore, it follows that the number of configurations for $(\cR,\cB)$ that yield $Z_{1,N}(\cR,\cB)=2p$ is $2\binom{N-1}{p}\binom{N-1}{p-1}$.\label{case:2}

\begin{figure}[t] 
    \centering
  \begin{minipage}[t]{\linewidth}
  \centering
    \includegraphics[width=0.48\linewidth]{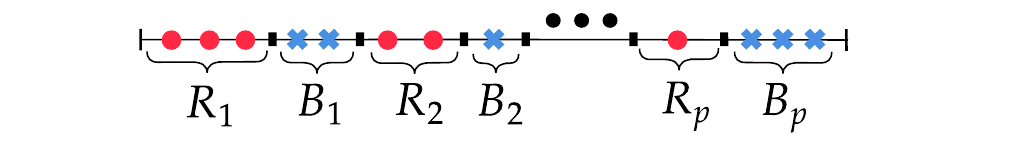}
    \includegraphics[width=0.48\linewidth]{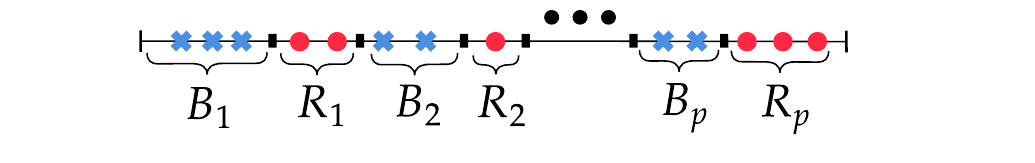} 
   \caption{Examples for case 1 in the proof of \cref{prop:3}.}
   \label{fig:k12}
  \end{minipage} 
  \vspace{0.5cm}
  
  \begin{minipage}[t]{\linewidth}
    \centering
    \includegraphics[width=0.48\linewidth]{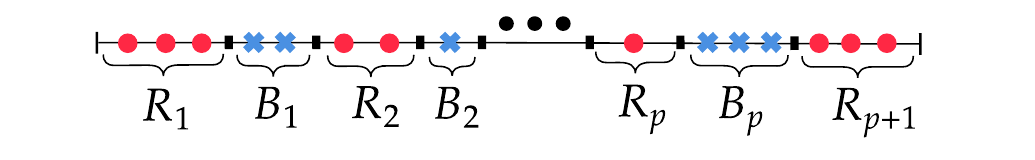}
    \includegraphics[width=0.48\linewidth]{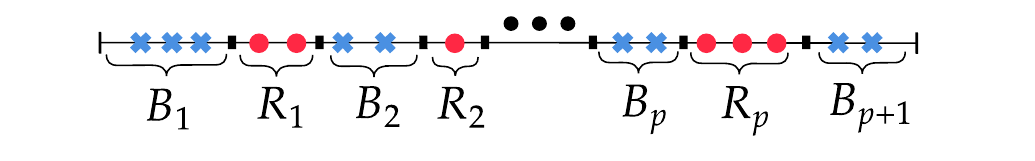}
   \caption{Examples for case 2 in the proof of \cref{prop:3}.}
   \label{fig:k34}
  \end{minipage} 
\end{figure}
Finally, the total number of possibilities is exactly the number of ways to choose $N$ points out of $2N$, which is given by the central binomial coefficient $\binom{2N}{N}$. Indeed, 
\begin{equation*}
   2\sum_{p=1}^N \binom{N-1}{p-1}^2   +2\sum_{p=1}^{N-1} \binom{N-1}{p} \binom{N-1}{p-1}= \binom{2N}{N}.
\end{equation*}
\end{proof}

\begin{remark}\label{rem:hypgeom}
   The probability mass function in \cref{eq:distrib1d} can be expressed in terms of the hypergeometric 
 distribution, whose mass function is defined by
\begin{equation*}
H(k;M,K,n)=\frac{\binom{K}{k}\binom{M-K}{n-k}}{\binom{M}{n}},
\end{equation*}
for all $0\leq M$, $0\leq K,n\leq M$, and $\max\{0,-M+K+n\}\leq k\leq \min\{K,n\}.$ Thus, for each $1\leq k\leq 2N-1$ we deduce
\begin{align*}
\mathbb{P}(Z_{1,N}=k)=
\begin{cases}
\displaystyle
N(2N-1)^{-1}H(p-1;2N-2,N-1,N-1), & \text{if } k=2p-1,\vspace{3mm}\\
\displaystyle
(N-1)(2N-1)^{-1}H(p;2N-2,N-1,N), & \text{if } k=2p.
\end{cases}
\end{align*}
\end{remark}

As a consequence of \cref{prop:3}, we get the distribution of $Z^\perp_{d,N}$ for any $d\geq1$:

\begin{corollary}\label{cor:final}Let $Z^\perp_{d,N}$ be as defined in \cref{eq:ZdNorth}, for some $d,N\geq1$.  For $1\leq k\leq2N-1$, we have
\begin{equation}\label{eq:Zgeqk}
\mathbb{P}(Z^\perp_{d,N}\geq k)
=\left(\sum_{p=\lceil\frac{k+1}{2}\rceil}^N\binom{N-1}{p-1}^2+\sum_{p=\lceil\frac{k}{2}\rceil}^{N-1}\binom{N-1}{p}\binom{N-1}{p-1}\right)^d2^d\binom{2N}{N}^{-d}.
\end{equation}
\end{corollary}

\begin{proof}
Let $k\in\{1,\dots,2N-1\}$. By definition of $Z^\perp_{d,N}$, and the fact that $Z^i_{d,N}$ are independent and identically distributed to $Z_{1,N}$ for all $i$, we can compute 
\begin{eqnarray*}
\mathbb{P}\left(Z^\perp_{d,N}\geq k\right)=\mathbb{P}\left(\min\limits_{i=1,\dots,d}Z_{d,N}^i\geq k\right)=\left(\mathbb{P}\left(Z_{1,N}\geq k\right)\right)^d.  
\end{eqnarray*}
We conclude the proof by applying \cref{prop:3}  to deduce
\begin{equation*}
\mathbb{P}\left(Z_{1,N}\ge k\right)
=\left(\sum_{p=\lceil\frac{k+1}{2}\rceil}^N\binom{N-1}{p-1}^2+\sum_{p=\lceil\frac{k}{2}\rceil}^{N-1}\binom{N-1}{p}\binom{N-1}{p-1}\right)2\binom{2N}{N}^{-1}.\end{equation*} 
\end{proof}

\smallskip

\begin{figure}[t]  
\centering
\includegraphics[width=0.35\linewidth]{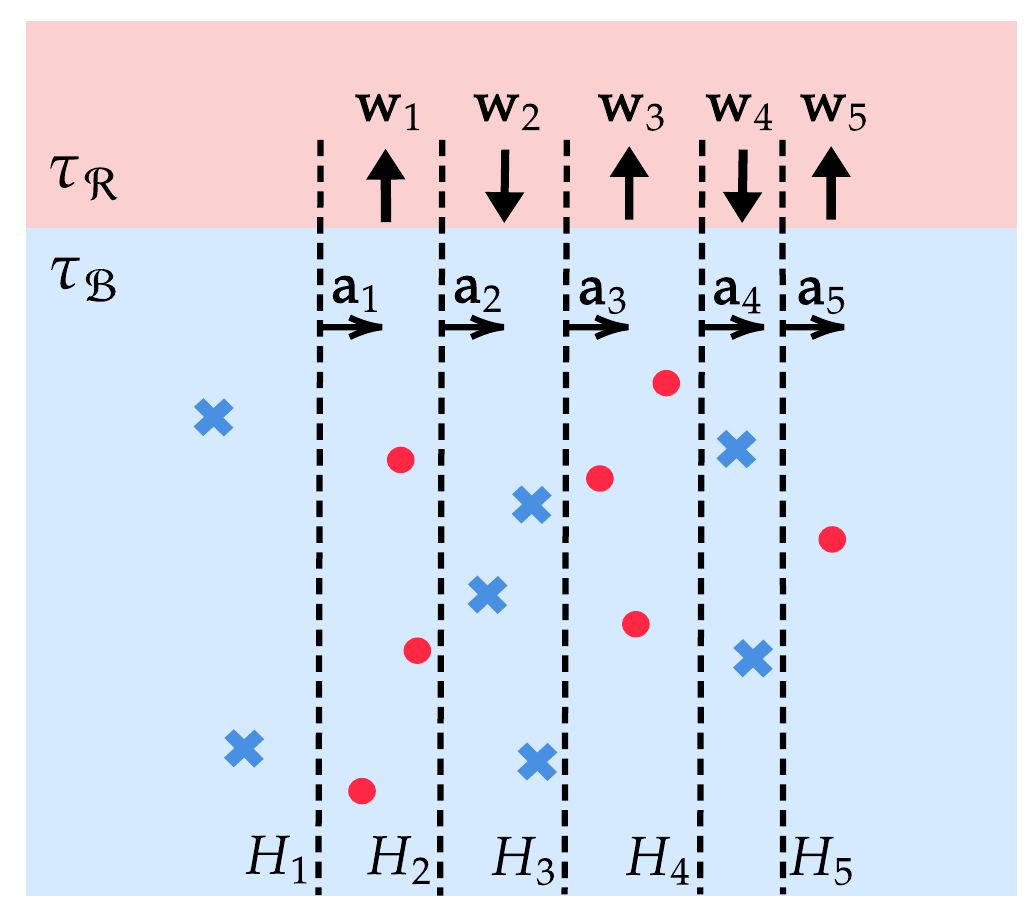}
\caption{Schematic overview of the full control method used to prove  \cref{thm2}. Parallel hyperplanes $H_i$ separate the two classes into point clusters, which then move vertically in  zigzag until all points are classified.}
\label{fig:algo1}
\end{figure}

The results derived in this section are now used to prove \cref{thm2}. For clarity, the whole control method is represented in \cref{fig:algo1} and formalized in \cref{alg:thm2}.

\begin{proof}[Proof of \cref{thm2}]
    By definition of \cref{eq:ZdNorth}, there exists a family of  hyperplanes \begin{equation*}
   H_1,\dots, H_z,  \hspace{1cm}\text{with }z=Z_{d,N}^\perp(\cR,\cB),  
    \end{equation*}
  that are orthogonal to $\bfe_i$ for some $i\in\{1,\dots,d\}$ and separate $(\cR,\cB)$. Moreover, $Z_{d,N}^\perp(\cR,\cB)$ follows the distribution given by \cref{eq:Zgeqk}. 
    
    Without loss of generality we can assume $i=1$, so the family of hyperplanes is
    \begin{equation*}\left\{H_j\;:\; x^{(1)}=h_j\right\}_{j=1}^{z} \quad \text{for some }0< h_1<\cdots<h_{z}<1.\end{equation*}
    For any fixed $k\in\{2,\dots,d\}$, we define $\tau_\cR=\{x^{(k)}>1\}$ and  $\tau_\cB= \{x^{(k)}\leq1\}$. We then classify the points by clusters based on their $x^{(1)}$-coordinates in ascending order. We can assume that $\{\bfx\in \cR\cup \cB:x^{(1)}<h_1\}\subset \cB$. If this is not the case, we swap the definitions of $\tau_\cR$ and $\tau_\cB$, and also the roles of $\cR$ and $\cB$ in this proof.

Taking $t_0=0$, we build the controls \begin{align}\label{eq:controls}
\notag&(\bfw,\bfa,b)(t)=\sum\limits_{j=1}^{z}(\bfw_j,\bfa_j,b_j)\,\mathbbm{1}_{[t_{j-1},t_j)}(t),\\
&\text{with}\quad\bfw_{j}=(-1)^{j+1}\bfe_k,\; \bfa_{j}=\bfe_1,\; b_{j}=-h_{j}. 
    \end{align}
  Each time horizon $t_j\geq t_{j-1}$ is chosen so that every point $\bfx\in \cR\cup \cB$ with $h_{j}<x^{(1)}<h_{j+1}$ is mapped via $(\bfw_j,\bfa_j,b_j)$ to its corresponding target region $\tau_\cR$ or $\tau_\cB$. This is possible because the dataset is finite. Moreover, these movements do not affect points that have already been classified, since those satisfy $x^{(1)}<h_j$.

The described method classifies all points according to their labels, and the number of switches is $L=z-1$. Consequently, $L$ inherits the probability distribution of $Z_{d,N}^\perp$.
\end{proof}

\begin{algorithm}[t]
\caption{Classification of two $N$-point sets}
\label{alg:thm2}
\begin{algorithmic}[1]
\Require Two $N$-point sets $\cR, \cB \subset [0,1]^d$
\Ensure Classification of $\cR$ and $\cB$

\State $ \displaystyle u \gets e_i$
    \Comment{Optimal direction for canonical separability}
\State $\displaystyle \tau_\cR \gets \{x^{(k)} > 1\},\quad\tau_\cB \gets \{x^{(k)} < 1\}$ 
    \Comment{For any coordinate $k \neq i$}
\State $\displaystyle \cR \cup \cB \gets \text{Reorder}(\cR \cup \cB;\, i)$ 
    \Comment{Sort all points by ascending $i$-th coordinate}
\State $\cS\gets\{x_1\},\quad \text{anchor}\gets 0$
\For{$\displaystyle j \in \{2,\dots,2N\}$}
     \If{$y_{j-1}=y_{j}$}
    \State $\cS\gets \cS\cup\{x_j\}$
    \ElsIf{$x_{j-1} \in \cB$ \textbf{ and } $x_j \in \cR$}
        \While{$\cS \subset \tau_\cR$}
            \State $\displaystyle \cR \cup \cB \gets \text{neuralODE}_{T=1}\Bigl(\cR \cup \cB;\, w = -e_k,\, a=u,\, b = \text{anchor}\Bigr)$
        \EndWhile
                    \State $\cS\gets\{x_j\},\quad \text{anchor}\gets-u \cdot x_{j-1}$
    \ElsIf{$x_{j-1} \in \cR$ \textbf{ and } $x_j \in \cB$}
        \While{$\cS \subset \tau_\cB$}
            \State $\displaystyle \cR \cup \cB \gets \text{neuralODE}_{T=1}\Bigl(\cR \cup \cB;\, w = e_k,\, a=u,\, b = \text{anchor}\Bigr)$
        \EndWhile
            \State $\cS\gets\{x_j\},\quad \text{anchor}\gets-u \cdot x_{j-1}$
    \EndIf
\EndFor

\end{algorithmic}
\end{algorithm}

\begin{remark}\label{rem:bv} \cref{thm2} can be interpreted in terms of the total variation seminorm\begin{equation*}
    |\theta|_{\mathrm{TV}(0,T)} = \sup_P\sum_{i=1}^{|P|}\|\theta(t_i)-\theta(t_{i-1})\|,
\end{equation*}
where the supremum is taken over all finite partitions $P=\{0=t_0<t_1<\cdots<t_{|P|}=T\}$. In particular, for piecewise constant functions we have that $|\cdot|_{\mathrm{TV}(0,T)}$ equals the sum of the magnitudes of its jump discontinuities. Thus, the control $\theta$ given by \cref{eq:controls} will satisfy
\begin{equation*}|\theta|_{\mathrm{TV}(0,T)}\leq L\cdot\max_{1\leq j\leq L}\sqrt{\|\bfw_{j+1}-\bfw_{j}\|^2+|b_{j+1}-b_{j}|^2}\leq 
\sqrt{5}L.\end{equation*} 
From here, we can estimate the probability distribution of $|\theta|_{\mathrm{TV}(0,T)}$ via 
\begin{equation*}\mathbb{P}\left(|\theta|_{\mathrm{TV}(0,T)}\geq \lambda\right)\leq \mathbb{P}\left(L\geq \frac{\lambda}{\sqrt{5}}\right),\hspace{1cm} \text{for all }\lambda>0\end{equation*}
and then apply \cref{eq:distribL}. \end{remark}

\begin{remark}\label{rem:ext1}
If we allow $\bfa$ to be any vector in $\R^d$, then the constraint \cref{eq:condcansep} can be removed in \cref{thm2}, ensuring the existence of  $T$, $(\tau_{\cR},\tau_{\cB})$ and $\theta$ in every case, rather than almost surely. To accomplish this, we choose a new vector basis in which \cref{eq:condcansep} holds, and note that the conclusion of \cref{prop:3} remains valid in the new coordinates.


\end{remark}

\section{Classification by separability in general position}\label{sec:4}

Let $(\cR,\cB)$ be as in \cref{eq:RB}, and  assume for now that $|\cR|=|\cB|=N$. We recall \cref{def:gp} and the random variable $Z_{d,N}$, defined in \cref{eq:ZdN}, to introduce the following quantity:
\begin{equation*}
M(d,N)\coloneqq\max\left\{Z_{d,N}(\cR,\cB):\,(\cR,\cB) \text{ as \cref{eq:RB} in general position, $|\cR|=|\cB|=N$}\right\}.
\end{equation*}
Observe that for $d=1$, any set of points is naturally in general position, so the result from \cref{thm3} applies directly. Consequently, we have $M(1, N) = 2N - 1$ for all $N$, as stated in \cref{eq:maxnumd}. However, when $d>1$, the situation changes significantly. Assuming general position eliminates the pathological configurations described in \cref{thm3}, thereby reducing the maximum possible value of $Z_{d,N}$. 

\begin{figure}[t]
    \centering
    \begin{minipage}[b]{0.49\textwidth}
        \centering
        \includegraphics[width=0.48\linewidth]{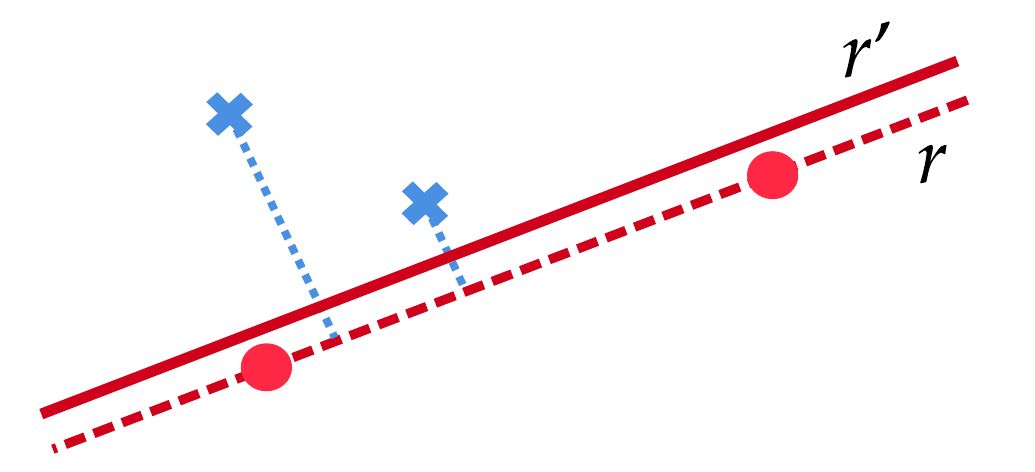}\;
        \includegraphics[width=0.48\linewidth]{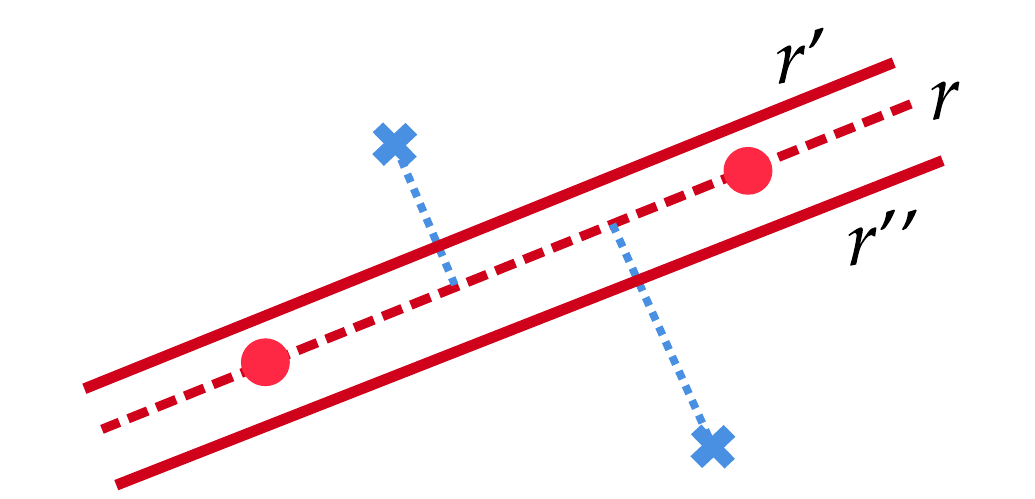}
        \caption{Separation of $(\cR,\cB)$ in general position in $\R^2$ with $|\cR|=|\cB|=2$ using at most two lines $r'$, $r''$.}
        \label{fig:gp1gp2}
    \end{minipage}
    \hfill
        \begin{minipage}[b]{0.49\textwidth}
        \centering
        \includegraphics[width=0.8\linewidth]{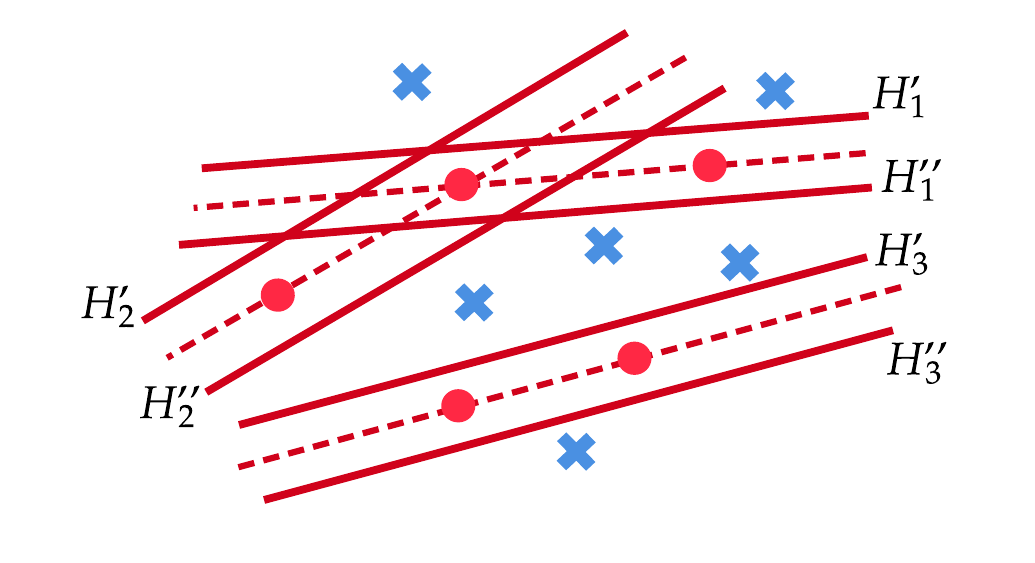}
    \caption{Hyperplanes constructed following the proof of \cref{thm:1gp}, separating two unbalanced classes.}
    \label{fig:exmethod}
    \end{minipage}    
\end{figure}

For example, assume that $d = N = 2$. Here, we find $M(2, 2) = 2 $, which is less than $2 \cdot 2 - 1 = 3$. To show this, consider connecting the two red points with a line $r$. Since all points are in general position, none of the blue points can lie on $r$. There are two possible cases: 
\begin{enumerate}
    \item If the blue points are on the same side of $r$ then we can separate $\cR$ from $\cB$ with one line $r'$ parallel to $r$, as shown in \cref{fig:gp1gp2} (left).

    \item If the blue points are on different sides of $r$ then we can separate $\cR$ from $\cB$ with two lines $r'$ and $r''$ parallel to $r$, as shown in \cref{fig:gp1gp2} (right).
\end{enumerate}
In the example of \cref{fig:gp1gp2} (left), 
 if the lines were restricted to be perpendicular to the canonical axes, two lines would be needed to separate $(\cR,\cB)$. Moreover, if the lines also had to be parallel---as in the framework of \Cref{sec:3}---then three lines would be required.

The argument used in the simple case $(d,N)=(2,2)$ can be extended to derive a general bound for $M(d,N)$ that shows the improvement over \cref{thm3} when the points are in general position. To this end, we cover $\cR$ with (possibly overlapping) subsets of size $d$, and then separate $(\cR,\cB)$ by isolating these subsets using hyperplanes under a transversality condition.

\begin{proposition}\label{thm:1gp} 
Let $d,N\geq1$, fix $i\in\{1,\dots,d\}$ and let $(\cR,\cB)$ be as \cref{eq:RB} in general position, with $|\cR|=|\cB|=N$. Then, there exist  hyperplanes $H_1',H_1'',\dots,H_{\lceil N/d\rceil}',H_{\lceil N/d\rceil}''$ that separate~$(\cR,\cB)$. Moreover, for all $j=1,\dots,\lceil N/d\rceil$ the following hold:
\begin{enumerate}
    \item $H_j'$ and $H_j''$ are parallel;
    \item $H_j'$ and $H_j''$ enclose exactly $\min\{d,|\cR|\}$ points of $\cR$ and no points of $\cB$.
    \item $H_j'$ and $H_j''$ are not orthogonal to $\bfe_i$;
\end{enumerate}
In particular, for all $d,N\geq1$ it follows that
\begin{equation}\label{eq:maxboundgenpos}
M(d,N)\leq 2\left\lceil \frac{N}{d}\right\rceil.
\end{equation}
\end{proposition}
\begin{proof} 

If $d\leq N$, we can choose $\hat R_1,\dots,\hat R_{\lceil N/d\rceil}\subset\cR$ such that $\cR=\hat R_1\cup\dots\cup \hat R_{\lceil N/d\rceil}$, with $|\hat R_j|=d$ for all $j$. Otherwise, build  $\hat R_1$ by augmenting the set $\cR$ with $d-N$ points chosen from $\R^d\setminus \cB$ so that $\hat R_1 \cup \cB$ remains in general position. 
Since $\cR\cup\cB$ is in general position, each $\hat R_j$ spans a unique hyperplane $\hat H_j\subset\R^d$ satisfying 
\begin{equation}\label{eq:spanhyp}
   (\cR\cup\cB\setminus \hat R_j) \cap \hat H_j = \emptyset\hspace{1cm}\text{and}\hspace{1cm} \hat H_j\neq \hat H_k\quad\text{ for all }j\neq k.
\end{equation} 
For each $j$, choose two hyperplanes $\hat H_j'$ and $\hat H_j''$ that are parallel to $\hat H_j$ and lie in the two distinct connected components of $\R^d\setminus \hat H_j$. If these hyperplanes are chosen sufficiently close to $\hat H_j$, then the region between $\hat H_j'$ and $\hat H_j''$ contains the subset $\hat R_j$ and no other points of $\cR\cup\cB$.  Consequently, $(\cR,\cB)$ is separated by the family of pairwise parallel hyperplanes
\begin{equation*}
\hat H_1',\hat H_1'',\dots,\hat H_{\lceil N/d\rceil}',\hat H_{\lceil N/d\rceil}''\subset\R^d.
\end{equation*}
By construction, these hyperplanes meet the first two conditions of the statement; however, they may not satisfy the third condition, as some might be orthogonal to $\bfe_i$. Suppose that exactly $2p$ hyperplanes are orthogonal to $\bfe_i$ for some $1 \le p \le \lceil N/d \rceil$. Without loss of generality, assume these are \begin{equation*}
\hat H_1',\hat H_1'',\dots,\hat H_p',\hat H_p''.
\end{equation*}
Since these hyperplanes are parallel, and different by \cref{eq:spanhyp}, we get $\hat R_i\bigcap \hat R_j=\emptyset$ for $i\neq j$.

Let $1\leq j\leq p$. Because $\cR \cup \cB$ is finite, we can slightly adjust the direction vector of $\hat H_j'$ and $\hat H_j''$. This yields new parallel hyperplanes $H_j'$ and $H_j''$ in $\R^d$ that are no longer orthogonal to $\bfe_i$, yet still enclose exactly $R_j$, with no other points of $\cR\cup\cB$ in between. For $p<j\leq\lceil N/d\rceil$, set $H_j'=\hat H_j'$ and $H_j''=\hat H_j''$. This yields a new family
\begin{equation*}
H_1', H_1'',\dots,  H_p', H_p'', H_{p+1}',H_{p+1}'',\dots,H_{\lceil N/d\rceil}', H_{\lceil N/d\rceil}''\subset\R^d
\end{equation*}
that meets all the required conditions. Moreover, this proves \cref{eq:maxboundgenpos}.


\end{proof}

\begin{figure}[t]  
\centering
\includegraphics[width=0.35\linewidth]{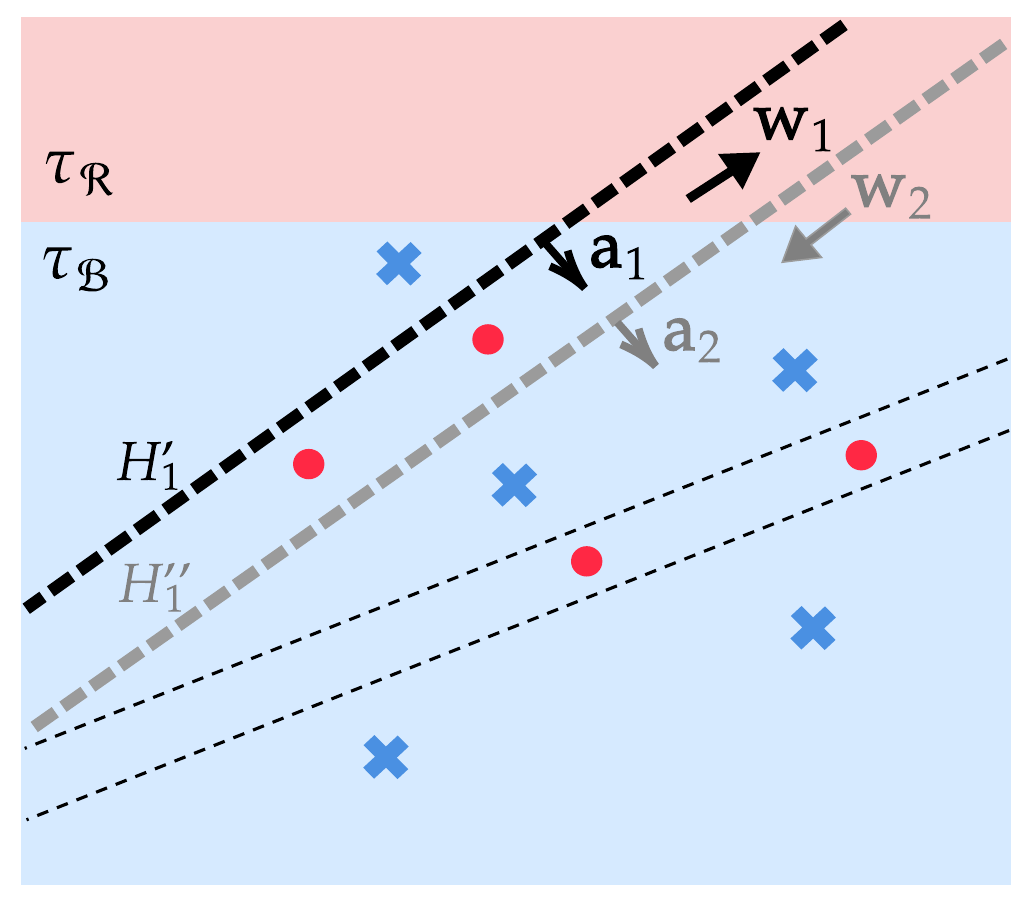}
\caption{Schematic overview of an iteration in the control method used to prove \cref{thm:trunnodes}. We aim to move the two red points that lie between the hyperplanes $H_1'$ and $H_1''$ to the region $\tau_{\cR}$ through a two-step process. First, we move all points within the half-space $H_1'^{+}$, determined by $\bfa_1$, in the direction of $\bfw_1$. Second, we move all points within the half-space $H_1''^{+}$, determined by $\bfa_2$, in the direction of $\bfw_1$. The controls are adjusted to ensure that the vector field in $H_1''^{+}$ maintains a unit norm at all times. Consequently, the points in $H_1''^{+}$ return exactly to their original positions after these two steps.}
\label{fig:algo2}
\end{figure}

 In the unbalanced case, where $|\cR| \neq |\cB|$, the method can be similarly applied to isolate the smaller set using pairwise parallel hyperplanes, as shown in \cref{fig:exmethod}. It follows:
\begin{corollary}\label{cor:unbalanced}
Let $d\geq1$. For any dataset $(\cR,\cB)$ as in \cref{eq:RB} in general position, there exists a family of $2\left\lceil\frac{ \min\{|\cR|,|\cB|\}}{d}\right\rceil$ hyperplanes that separate the points by labels and satisfy the three conditions of \cref{thm:1gp}, with $\cR$ replaced by $\operatorname{argmin}_{\cR,\cB}\left\{|\cR|,|\cB|\right\}$.
\end{corollary}
\Cref{cor:unbalanced} enables the separation of any finite dataset in general position in $\R^d$ into clusters of size $d$. We now describe an inductive approach to classify these clusters using the dynamics of neural ODEs. First, let us consider the system
\begin{equation}\label{eq:nodetrun}
             \dot \bfx (t)   = \bfw(t)\,\sigma_{\text{\rm trun}} (\bfa(t)\cdot\bfx(t)+b(t)),
\end{equation}
where the activation function is defined as
\begin{equation}\label{eq:trunrelu}
\sigma_{\text{\rm trun}}(z)\coloneqq\min\big\{1,(z)_+\big\} = (z)_+ - (z-1)_+,\hspace{1cm}\text{for }z\in\R.
\end{equation}
Observe that the flow map of \cref{eq:nodetrun} remains well-defined because $\sigma_{trun}$ is Lipschitz-continuous. The introduction of $\sigma_{trun}$ is motivated by its properties, which facilitate the inductive argument. Specifically, in the half-space defined by $\bfa(t)\cdot\bfx+b(t)>1$, the system \cref{eq:nodetrun} simplifies to $\dot\bfx = \bfw(t)$. This property can be used to ensure that points already classified will remain so after each inductive step. For more details, see \cref{fig:algo2}, which serves as support for the proof, and the formalized control method in \cref{alg:thm1}.

\begin{proposition}\label{thm:trunnodes}
Let $d\geq 2$. For any dataset $(\cR,\cB)$ defined as in \cref{eq:RB} in general position, and any pair of target sets $(\tau_\cR,\tau_\cB)$ defined as in \cref{eq:targets}, there exist $T>0$ and a piecewise constant control 
$\theta\in\Theta_T$ whose number of discontinuities is  \begin{equation*}
L=2\left\lceil \frac{\min\left\{|\cR|,|\cB|\right\}}{d}\right\rceil-1,    
\end{equation*}
 such that the flow map of \cref{eq:nodetrun} satisfies    $\Phi_T(\cR;\theta)\subset \tau_\cR$ and $\Phi_T(\cB;\theta)\subset \tau_\cB$.  
\end{proposition}
\begin{proof}
Let $i\in\{1,\dots,d\}$ be fixed and consider $\tau_\cR=\{x^{(i)}>1\}$,  $\tau_\cB=\{x^{(i)}\leq 1\}$. The strategy is to evolve $\cR$ into the interior of $\tau_\cR$ while keeping fixed $\cB$. We first assume $|\cR|=|\cB|=N$ and will then extend to the general case $|\cR|\neq|\cB|$. 

By \cref{thm:1gp}, there exist $2\lceil N/d\rceil$ pairwise parallel hyperplanes separating $(\cR,\cB)$. Assume these hyperplanes are defined by \begin{equation*}
H_j':\bfu_j \cdot \bfx + h_j' = 0,\qquad\text{and}\qquad H_j'':\bfu_j \cdot \bfx + h_j'' = 0, \hspace{1cm} \text{for }j = 1,\dots, \left\lceil \frac{N}{d} \right\rceil,  
\end{equation*}
where $h_j'>h_j''$, and $\bfu_j\in\mathbb{S}^{d-1}$ satisfies $|\bfu_j\cdot\bfe_i|<1$.  Moreover, the region between each pair $(H_j',H_j'')$ encloses a subset $R_j \subset \cR$ such that $|R_j|=\min\{d,|\cR|\}$, and no other points of $\cR \cup \cB$ lie between these hyperplanes, namely,
\begin{equation*}
(H_j'^+\setminus H_j''^+)\cap \left(\cR\cup\cB\right) = R_j\subset\cR\hspace{1cm}\text{with }|R_j|=\min\{d,|\cR|\},
\end{equation*}
where $H_j'^+=\{\bfx\in\R^d:\bfu_j\cdot \bfx+h_j'>0\}$ and $H_j''^+=\{\bfx\in\R^d:\bfu_j\cdot \bfx+h_j''>0$\}. Note that $H_j''^+\subset H_j'^+$ because $h_j'>h_j''$. Taking $t_0=0$, we define \begin{equation*}(\bfw,\bfa,b)=\sum\limits_{k=1}^{2\lceil N/d\rceil}(\bfw_k,\bfa_k,b_k)\,\mathbbm{1}_{(t_{k-1},t_k)}\end{equation*}
such that, for $j=1,\dots,\lceil N/d\rceil$:
\begin{enumerate}
    \item $(\bfw_{2j-1},\bfa_{2j-1},b_{2j-1})=(\bfv_j,\bfu_j/d_j',h_j'/d_j')$, where \begin{itemize}
        \item $\bfv_j\in\mathbb{S}^{d-1}$  satisfies $\bfv_j\cdot\bfu_j=0$ and $\bfv_j\cdot\bfe_i>0$ (for instance, take $\bfv_j$ to be the normalized projection of $\bfe_i$ onto the orthogonal subspace $\langle \bfu_j \rangle^{\perp}$);
        \item $d_j'=\min\left\{\sigma_{\text{\rm trun}}(\bfu_j\cdot \bfx+h_j'):\bfx\in R_j\right\}>0$.
\end{itemize}
Inside the half-space $\{\bfx\in\R^d:\bfu_j\cdot\bfx+h_j'\geq d_j'\}\subset H_j'^+$ and over $(t_{2j-2},t_{2j-1})$, equation \cref{eq:nodetrun} becomes $\dot\bfx=\bfv_j$ with $\bfv_j\cdot\bfe_i>0$. Since the dataset is finite and contained in this half-space, we can choose $t_{2j-1}>t_{2j-2}$ such that $\Phi_{t_{2j-1}}(\bfx)\in \tau_\cR$ for all $\bfx\in R_j$.

\item $(\bfw_{2j},\bfa_{2j},b_{2j})=(-\bfv_j,\bfu_j/d_j'',h_j''/d_j'')$, where \begin{itemize}
        \item $d_j''=\min\left\{\sigma_{\text{\rm trun}}(\bfu_j\cdot \bfx+h_j''):\bfx\in (\cR\cup\cB)\cap H''^+_j\right\}>0.$
    \end{itemize}
    Inside the half-space $\{\bfx\in\R^d:\bfu_j\cdot\bfx+h_j''\geq d_j''\}\subset H_j''^+$ and over $(t_{2j-1},t_{2j})$, equation \cref{eq:nodetrun} becomes $\dot\bfx=-\bfv_j$. Now, we set $t_{2j}=2t_{2j-1}-t_{2j-2}$ so that \begin{equation*}
   \Phi_{t_{2j}}(\bfx) =\Phi_{t_{2j-2}}(\bfx)\hspace{1cm}\text{for all }\bfx\in (\cR\cup \cB) \cap H_j''^+.    
    \end{equation*}
    All the while, we have 
    \begin{equation*}
   \Phi_{t_{2j}}(\bfx) = \Phi_{t_{2j-1}}(\bfx)\hspace{1cm}\text{for all }\bfx\in (\cR\cup \cB) \cap H_j''^-.    
    \end{equation*}
 
\end{enumerate}
For $T=t_{2\lceil N/d\rceil}$, we conclude 
$\Phi_T(\cR)\subset~\tau_\cR$ and $\Phi_T(\cB)=\cB\subset \tau_\cB$, with $L=2\lceil N/d\rceil-1$.

If $|\cR|<|\cB|$, we apply the same argument to obtain $\Phi_T(\cR)\subset~\tau_\cR$  and $\Phi_T(\cB)=\cB\subset \tau_\cB$ with $L=2\lceil |\cR|/d\rceil-1$, by virtue of \Cref{cor:unbalanced}. If $|\cR|>|\cB|$, we swap the roles of $\cR$ and $\cB$ and the definitions of $\tau_\cR$ and $\tau_\cB$, with $L=2\lceil |\cB|/d\rceil-1$.
\end{proof}

The following lemma formalizes an idea from \cite[section 2.2]{domzuazuaNormFlows}. It shows that any flow of \cref{eq:nodetrun} can be represented as the composition of two flows of \cref{eq:node}, if $\bfw$ and $\bfa$ are orthogonal:
\begin{lemma}\label{lem:trunintorelu}
 Let $\theta=(\bfw,\bfa,b)\in\mathbb{S}^{d-1}\times\R^{d+1}\times\R$ with $\bfw\cdot\bfa=0$. Then, $\theta_1=(\bfw,\bfa,b)$ and $\theta_2 = (-\bfw,\bfa,b-1)$ satisfy 
    \begin{equation*}
        \Phi_{t}(\Phi_{t}(\cdot\,;\,\theta_1)\,;\,\theta_2)=\Phi_t^{\text{\rm trun}}(\cdot\,;\,\theta)\hspace{1cm}\text{for all }t>0,
    \end{equation*}
    where $\Phi_t$ and $\Phi_t^{\text{\rm trun}}$ are the flow maps of \cref{eq:node} and \cref{eq:nodetrun}, respectively.
\end{lemma}
\begin{proof}
Since $\bfw\cdot\bfa = 0$, it holds $$\frac{d}{dt}(\bfa\cdot\bfx(t)+b)=\bfa\cdot\dot\bfx(t)  = \bfa\cdot\bfw\sigma(\bfa\cdot\bfx(t)+b)= 0.$$
Thus, $\sigma(\bfa\cdot\bfx(t)+b)$ is constant for all $t$. Suppose $\bfa\cdot\bfx+b\leq1$. Then \begin{equation*}
    \Phi_t^{\text{\rm trun}}(\bfx\,;\,\theta) = \bfx+t\bfw(\bfa\cdot\bfx+b)=\Phi_{t}(\bfx\,;\,\theta_1)=\Phi_{t}(\Phi_{t}(\bfx\,;\,\theta_1)\,;\,\theta_2).   
    \end{equation*}
    Otherwise, if $\bfa\cdot\bfx+b>1$,
    \begin{align*}
    \Phi_{t}(\Phi_{t}(\bfx\,;\,\theta_1)\,;\,\theta_2)& = \Phi_{t}(\bfx+t\bfw(\bfa\cdot\bfx+b)\,;\,\theta_2)\\&=\bfx+t\bfw(\bfa\cdot\bfx+b)-t\bfw(\bfa\cdot(\bfx+t\bfw(\bfa\cdot\bfx+b))+b-1)\\
&=\bfx+t\bfw=\Phi_t^{\text{\rm trun}}(\bfx\,;\,\theta).
    \end{align*}
\end{proof}
We can now use \cref{lem:trunintorelu}, together with \cref{thm:trunnodes}, to demonstrate \cref{thm1}.

\begin{algorithm}[t]
\caption{Classification of two point sets via \cref{eq:nodetrun}}
\label{alg:thm1}
\begin{algorithmic}[1]
\Require Two finite sets $\cR,\cB\subset[0,1]^d$ in general position; Direction $e_i$
\Ensure Classification of $\cR$ and $\cB$

\State $\cS \gets \operatorname{argmin} \{|\cR|, |\cB|\},\quad \tau_\cS \gets\{x^{(i)}>1\}$  
\State $\mathcal{C} \gets \text{cover}(\cS, d, \lceil |\cS|/d \rceil)$  \Comment{Covering for $\cS$ by $\lceil |\cS|/d \rceil$ subsets $\cS_j$ with $|\cS_j|=d$}

\For{$(j,\cS_j)$ in $\text{enumerate}(\mathcal{C})$}
    \State $H_j(a_j, b_j) = \{x:a_j \cdot x + b_j = 0,\,\|a_j\|=1\} \gets \text{Span}(\cS_j)$
    \If{$|a_j\cdot e_i|=1$}
    \State $a_j \gets \text{NormalizedPerturbation}(a_j)$ \Comment{Achieves $|a_j\cdot e_i|<1,\,\|a_j\|=1$}
    \EndIf
    \State $\delta_j \gets \min\{\operatorname{dist}(x, H_j):\, x\in (\cR\cup\cB)\setminus\cS_j\}$  
    \State $b_j' \gets b_j + 0.5\,\delta_j,\quad b_j'' \gets b_j - 0.5\,\delta_j$ 
    \State $w_j \gets \text{NormalizedProjection}_{\langle a_j\rangle^\perp}(e_i)$
    \State $d_j' = \min_{x\in \cS_j}\sigma_{\text{trun}}(a_j\cdot x+b_j')$
       \State $d_j'' = \min_{x\in (\cR\cup\cB)\cap\{x\,:\,a_j\cdot x+b_j''>0\}}\sigma_{\text{trun}}(a_j\cdot x+b_j'')$
    \While{\(\cS_j\not\subset\tau_{\cS}\)}
        \State \(\cR\cup\cB \gets \text{neuralODE}_{T=1}(\cR\cup\cB\,; \,w=w_j,\,a=a_j/d_j',\,b=b_j'/d_j')\)
        \State \(\cR\cup\cB \gets \text{neuralODE}_{T=1}(\cR\cup\cB\,; \,w=-w_j,\,a=a_j/d_j'',\,b=b_j''/d_j'')\)
    \EndWhile
\EndFor
\end{algorithmic}
\end{algorithm}

\begin{proof}[Proof of \cref{thm1}]
 By virtue of \cref{thm:trunnodes}, there exist $\bar T>0$ and a piecewise constant control $\bar\theta=(\bar\bfw,\bar\bfa,\bar b)\in\Theta_{\bar T}$, which presents $2\lceil \min\{|\cR|,|\cB|\}/d\rceil-1$ discontinuities, such that \begin{equation*}
\Phi_{\bar T}^{\text{\rm trun}}(\cR;\bar\theta)\subset\tau_\cR\hspace{1cm}\text{and}\hspace{1cm}\Phi_{\bar T}^{\text{\rm trun}}(\cB;\bar\theta)\subset\tau_\cB, 
\end{equation*}
where $\Phi_t^{\text{\rm trun}}$ is the flow map of \cref{eq:nodetrun}. Moreover, $\bar\bfw(t)\cdot\bar\bfa(t)=0$ for all $t$, so by \cref{lem:trunintorelu} there exists a piecewise constant  $\theta\in\Theta_{2\bar T}$ with $4\lceil \min\{|\cR|,|\cB|\}/d\rceil-1$ discontinuities, such that $\Phi_{\bar T}^{\text{\rm trun}}(\cdot\,;\,\bar\theta)=\Phi_{2\bar T}(\cdot\,;\,\theta)$.  We conclude by defining $T=2\bar T$.
\end{proof}
\begin{remark}
Similarly to \cref{rem:bv}, we can extract from \cref{thm1} a bound for the total variation of the controls as
\begin{equation*}
    |\theta|_{\mathrm{TV}(0,T)}\leq 2L\|\theta\|_\infty\leq 2\sqrt{1+\|(\bfa,b)\|_\infty^2}\left(4\left\lceil \frac{\min\left\{|\cR|,|\cB|\right\}}{d}\right\rceil-1\right),
\end{equation*}
where $\|(\bfa,b)\|_\infty$ depends on the minimum separation between points of $\cR\cup\cB$ and the maximum distance of any point to the origin.
\end{remark}
The structure of pairwise parallel hyperplanes defined in \cref{thm:1gp} can be exploited for control when combined with certain architectures. We propose \emph{triangular neural ODEs}
\begin{equation}\label{eq:nodealt}
             \dot \bfx (t)  = 
             \bfw(t)\,\sigma_{\textrm{Tr}} (\bfa(t)\cdot\bfx(t)+b(t)),
\end{equation}
where the activation $\sigma_{\textrm{Tr}}$ is defined as 
\begin{equation}\label{eq:sigmatr}
\sigma_{\textrm{Tr}}(z)=\max\{ 1-|z-1| ,0\},\hspace{1cm}\text{for }z\in\R.
\end{equation}
 The function $\sigma_{\textrm{Tr}}$ remains globally Lipschitz, ensuring that the flow map is well-defined. 

By replacing $\sigma_{\text{\rm trun}}$ with  $\sigma_{\textrm{Tr}}$, we can improve upon \cref{thm:trunnodes}. Specifically, we can halve the number of values taken by the control $\theta$. Adapting the proof is straightforward: Step 2 can be omitted, as the support of $\sigma_{\textrm{Tr}}$ is restricted to the interval $(0,2)$.

\begin{corollary}\label{thm:nodealt}
Let $d\geq 2$. For any dataset $(\cR,\cB)$ defined as in \cref{eq:RB} in general position, and any pair of target sets $(\tau_\cR,\tau_\cB)$ defined as in  \cref{eq:targets}, there exist $T>0$ and a piecewise constant control
$\theta\in\Theta_T$ whose number of discontinuities is\begin{equation*}
    L=\left\lceil \frac{\min\left\{|\cR|,|\cB|\right\}}{d}\right\rceil-1,
\end{equation*}
 such that the flow map of the neural ODE \cref{eq:nodealt} 
satisfies      $\Phi_T(\cR;\theta)\subset \tau_\cR$ and $\Phi_T(\cB;\theta)\subset \tau_\cB$.
\end{corollary}

\begin{remark}
    Our separability-based methodology aligns with existing research in discrete geometry. In \cite{BoU1995}
    , the authors establish lower and upper bounds on the minimum number of hyperplanes required to individually separate $N$ points in general position in $\R^d$. In \cite{freimer_complexity_1991} it is shown that finding such number for an arbitrary point set is an NP-complete problem. From a computational perspective, efficient algorithms for separability in low-dimensional spaces are proposed in \cite{AHMSS,HOULE1993139}.
\end{remark}

\begin{remark}\label{rem:vc}
The VC dimension of a classification model is defined as the size of the largest set of points that can be shattered, meaning any arbitrary binary assignment of labels is possible \cite{vc}. We can interpret \cref{thm1} in terms of the VC dimension of our model.  Specifically, we can establish a lower bound: by fixing the maximum number of discontinuities in the controls to $L\geq1$, our method can shatter any set of $Ld$ points in $\R^d$.
\end{remark}


\section{Numerics}\label{sec:5}

We present a computational test\footnote{The code, based on PyTorch \cite{pytorch}, is available at \href{https://github.com/antonioalvarezl/2024-WCS-NODEs}{\texttt{https://github.com/antonioalvarezl/2024-WCS-NODEs}}} designed to evaluate the capacity of neural ODEs for binary classification. Specifically, our objective is to numerically estimate the minimal complexity they require to classify an arbitrary dataset of a fixed size, and compare it with the complexity we used in \cref{alg:thm2}.  

\begin{figure}[t]
\centering
\includegraphics[width=0.32\textwidth]{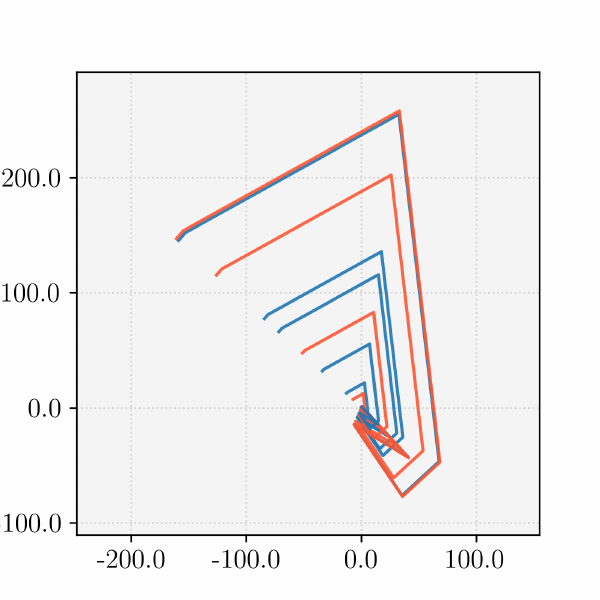}
\includegraphics[width=0.32\textwidth]{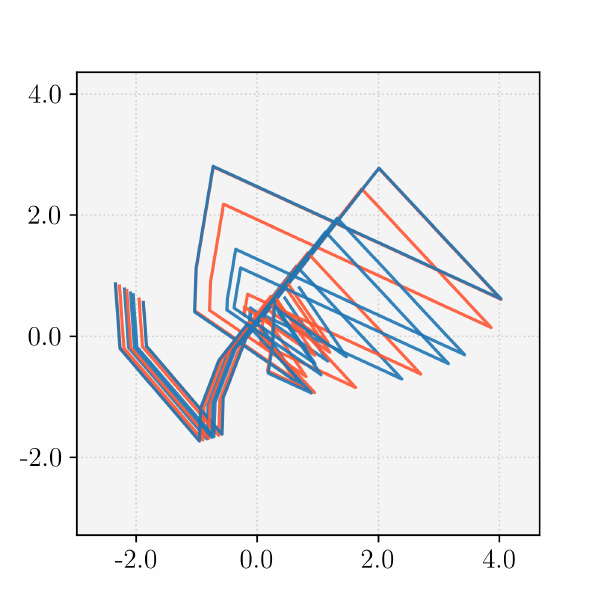}
\includegraphics[width=0.32\textwidth]{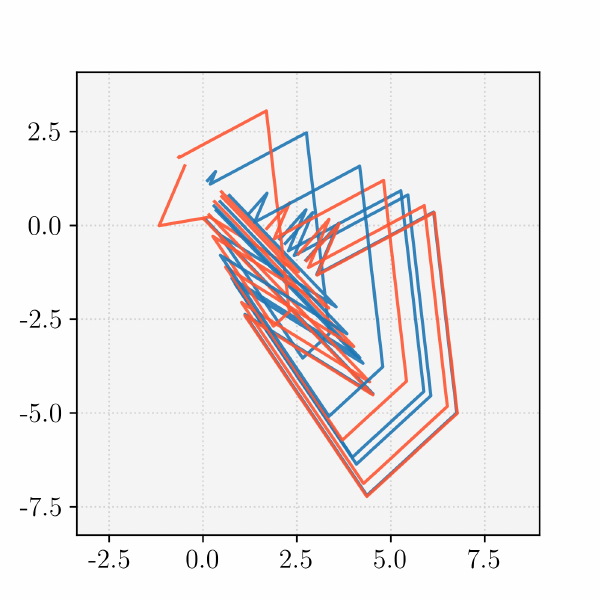}
    \caption{Trajectories for $\sigma=\operatorname{ReLU}$ exhibit an exponential drift when $\bfa\cdot\bfx+b>1$ (left). If $\sigma=\operatorname{tanh}$ (center) or $\sigma=\sigma_{\mathrm{trun}}$ (right), this drift is mitigated because $\|\sigma\|_\infty\leq 1$. Here, we have used $N=5$, $L=10$, and $T=60$.}
    \label{fig:trunvsrelu}
\end{figure}

Let $\scrD=\{(\bfx_n,y_n)\}\subset\R^d\times\{1,0\}$ be a finite dataset and $\cR=\{\bfx_n:(\bfx_n,1)\in\scrD\}$, $
 B=\{\bfx_n:(\bfx_n,0)\in\scrD\}$ be such that $|\cR|=|\cB|=N$. For a fixed $L\geq0$ and time $T>0$, we aim to find piecewise constant controls $\theta\in\Theta_T$ with $L$ discontinuities (or switches) such that the flow map of the neural ODE satisfies 
\begin{equation}\label{eq:obj}
\Phi_T(\cR;\theta)\subset \left\{x^{(1)}>1\right\}\hspace{1cm}\text{and}\hspace{1cm}\Phi_T(\cB;\theta)\subset \left\{x^{(1)}\leq1\right\}.
\end{equation}
We opt for $\sigma_{\text{\rm trun}}$ as in \cref{eq:trunrelu} and equation \cref{eq:nodetrun} to better visualize the results (see \cref{fig:trunvsrelu}). In this case, \cref{thm:trunnodes} establishes that $L=2\lceil N/d\rceil -1$ switches ensure any dataset in general position can be classified. This condition is  generically satisfied when all points are sampled from a non-singular probability measure, such as $U([0,1]^d)$.

We design our experimental setup in the following way:\newline
\textbf{Data}. 
We use ten random datasets, each consisting of $N = 30$ red points ($\cR$) and $N = 30$ blue points ($\cB$), sampled from $U([0,1]^d)$ for various dimensions $d$ ranging from 2 to $2N$.\newline
\textbf{Model}. We consider equation \cref{eq:nodetrun} with controls that are piecewise constant over $L+1$ time intervals of length $\Delta t=T(L+1)^{-1}$, where $T=60$. To approximate the solution, we use an explicit Euler discretization scheme with a step size of $0.25\Delta t$. The control values are initialized using Kaiming uniform initialization; that is, each component is sampled from $U([-\alpha,\alpha])$, where $\alpha$ is the inverse of the number of input units in the weight tensor.\newline
\textbf{Error}. 
To enforce the correct classification of all points as defined by \cref{eq:obj}, we introduce the following margin-based loss function:
\begin{equation*}\mathcal{L}(\theta)=\frac{1}{|\cR|}\sum_{\bfx_n\in \cR}\operatorname{dist}\left(\Phi_T(\bfx_n;\theta),\{x^{(1)}<1.5\}\right)^2 + \frac{1}{|\cB|}\sum_{\bfx_n\in \cB}\operatorname{dist}\left(\Phi_T(\bfx_n;\theta),\{x^{(1)}>0.5\}\right)^2
\end{equation*}
\textbf{Training.} Optimization is conducted using the Adam optimizer \cite{Kingma2014AdamAM} with a learning rate of 0.01, which was determined to be optimal through a grid search over a range of values.\newline
\textbf{Procedure}. First, for each $d$, we verify that setting $L=2\lceil N/d\rceil-1$ switches is sufficient to classify all ten datasets. Then, we gradually decrease the value of $L$ until at least one of the datasets fails to be classified, according to the stopping criteria defined below. \newline
\textbf{Stopping criteria}. If condition \cref{eq:obj} is met,  classification is successful and the training is stopped. Conversely, we consider the following three failure stopping criteria:
\begin{enumerate}
\item The maximum number of 70000 epochs is reached.
 \item Slow convergence, if 
 $\mathcal{L}\geq 0.15$ at 20000 epochs or $\mathcal{L}\geq0.1$ at 40000 epochs, or if the minimum error does not decrease over 5000 consecutive epochs.
 \item Local minima detection, if the maximum relative error over 50 consecutive epochs exceeds a threshold of $10^{-20}$.
\end{enumerate}
{\bf Close.} We conclude that the model with $L$ switches does not have the capacity to classify 60 points if any of the failure stopping criteria is met in 20 randomized initializations of the parameters for any of the ten datasets.

\medskip

 \begin{figure}[t!]  
\centering
\includegraphics[width=0.35\textwidth]{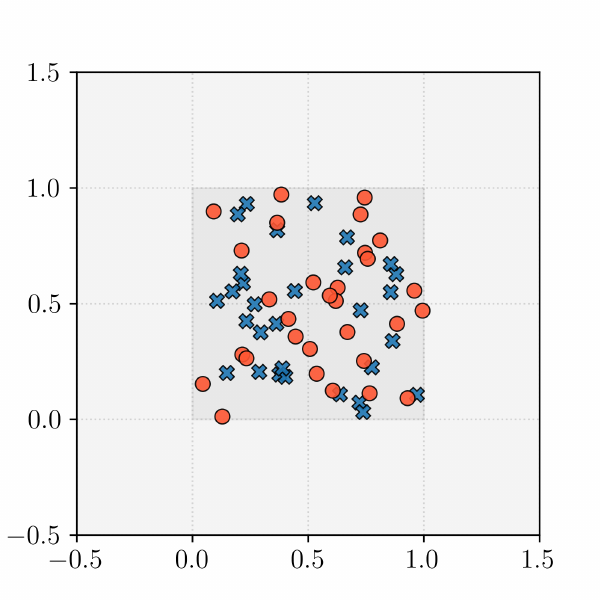}
\qquad\qquad
\includegraphics[width=0.35\textwidth]{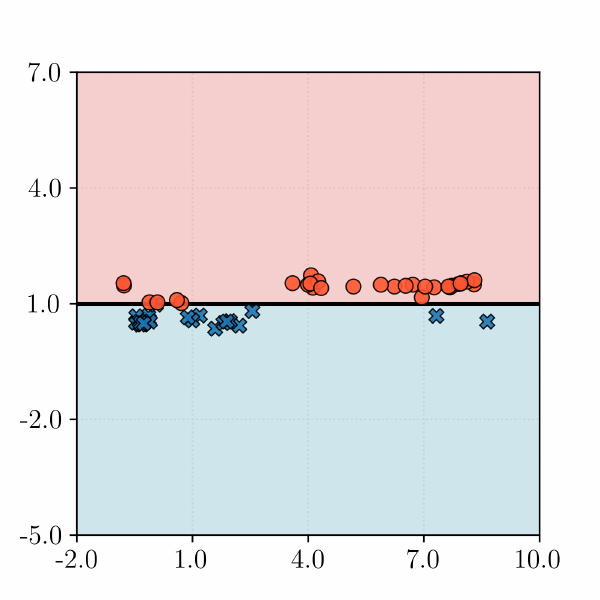}
\caption{Left: Initial data comprising $30$ red points and $30$ blue points in $[0,1]^2$. Right: Final positions at time $T=60$, having fixed $L=37$ switches.}
\label{fig:points}
\end{figure}

\begin{figure}[t]\centering\includegraphics[scale=0.4]{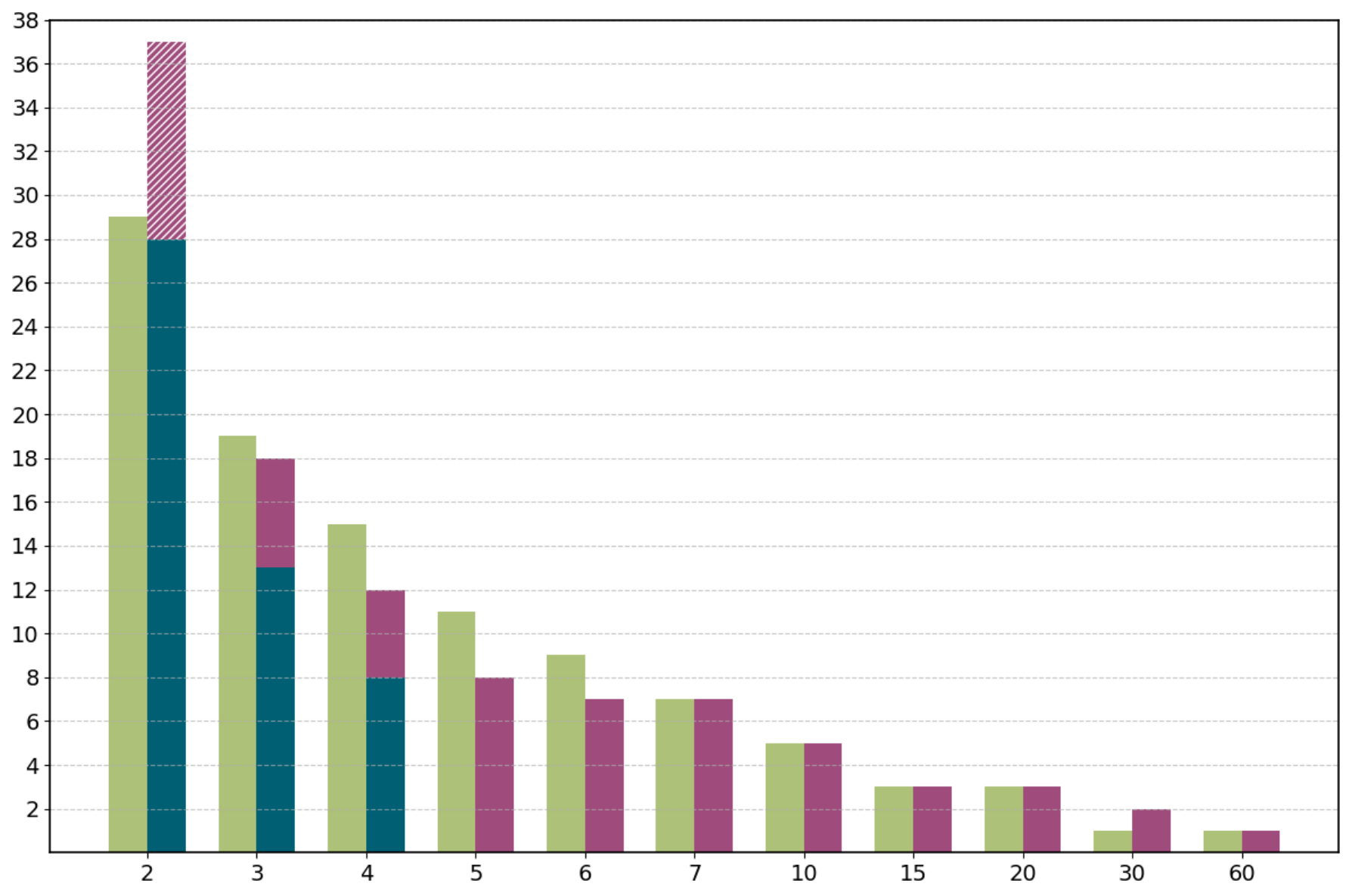} 
\caption{$L$ versus $d$ for fixed $N=30$. Green bars represent the number of switches in \cref{thm:trunnodes}, given by $2\lceil N/d\rceil -1$. 
Purple bars show the minimum number of switches required by gradient-based training to successfully classify all datasets. Dark bars indicate the minimum number of switches found after data rescaling.}    \label{fig:barplot}\end{figure}

The experiments were conducted using specific seed settings to ensure reproducibility, and the results are shown in \cref{fig:barplot}. We observe that $L=2\lceil N/d\rceil -1$ is typically close to the optimal value $L_*$ obtained through training. Specifically, for $d \in \{7, 10, 15, 20, 60\}$, we get $L_*=2\lceil N/d\rceil -1$. For $d\in\{2,30\}$, we even find that 
$L_*>2\lceil N/d\rceil -1$.

We observed that training requires an increasing number of parameter initializations in lower dimensions. This tendency is particularly pronounced when $d=2$, as shown by the purple bar with a hatched pattern. In that case, only five out of ten datasets were successfully classified using $L=37$ switches. We were unable to find any higher value of $L$ that could classify the remaining five datasets. 

Motivated by the difficulties found in low dimensions, we applied a preprocessing step to improve the separability of the data points. Specifically, for $d \in \{2, 3, 4\}$, we standardized the data by subtracting the mean and scaling to unit variance. This reduces the risk of the algorithm becoming trapped in poor local minima due to an overly dense point cloud. The results are represented by the blue bars of  \cref{fig:barplot}.


In summary, we observed that the value of 
$L=2\lceil N/d\rceil -1$ switches given in \cref{thm:trunnodes} is a suboptimal bound that closely approximates the optimal value obtained via gradient descent, especially for large $d$ with fixed $N$.

\section*{Conclusions}

In this work, we have explored the capacity of neural ODEs to classify an arbitrary dataset, using the framework of simultaneous control of clusters of points and geometric separability techniques to define these clusters.

In our main result, we have established a new bound of $L = 4 \left\lceil \min\{|\cR|, |\cB|\} / d \right\rceil $ control switches that guarantees classification of any dataset given by $ (\cR, \cB) $ in $\R^d$, using the single-neuron model \cref{eq:node}. The only assumption is that all points are in general position (see \cref{def:gp}), a condition that is satisfied almost surely by any finite set sampled from a non-singular measure. Our result improves the previous bound of $L = 3 \min\{|\cR|, |\cB|\} $, by leveraging the better separability that points present in higher dimensions. From a technical perspective, we demonstrate that any subset of $d$ points can be isolated between two parallel hyperplanes within any finite dataset in general position. We have complemented this result with some numerical experiments that supported the optimality of our bound, particularly in high dimensions. 

Our second result acknowledges that maximal complexity is rarely necessary, as data points of the same class are often initially close together, which facilitates their classification. We develop a new control method to derive the probability distribution of $L$ when all data points are i.i.d. and both classes have fixed size $N\geq1$. Our method once again emphasizes the advantage of high dimensionality. Specifically, we deduce that if $d \gtrsim 2^N / \sqrt{N}$ then classification using an autonomous neural ODE (characterized by $L=0$) is possible with high probability. Additionally, we have characterized the pathological configurations in which the maximum value of $L = 2N - 2$ occurs.

Our results can be seen as a manifestation of the \emph{blessing of dimensionality} \cite{bless}, a phenomenon where machine learning problems become more tractable in high dimensions. This contrasts with the well-known curse of dimensionality, which occurs when data sparsity in higher dimensions leads to an exponential increase in computational complexity and data requirements.

We now discuss some connections and extensions of our results:
\begin{itemize}
    \item \textbf{Multiclass classification}. For any $S \geq 2$, let $\cR_1,\dots, \cR_S \subset \R^d$ be finite subsets whose union is in general position. Suppose $|\cR_S|=\max|\cR_{i}|$. We first address binary classification between $\cR_1$ and $\cR_2 \cup \cdots \cup \cR_S$ using \cref{thm:trunnodes}. Once these two sets are linearly separable under the flow map $\Phi_T$, we proceed to classify $\Phi_T(\cR_2)$ against $\Phi_T(\cR_1 \cup \cR_3 \cup \cdots \cup \cR_S)$ in a similar manner. Since \cref{alg:thm1} consistently fixes the larger of the two sets, we can proceed inductively and ensuring that ultimately every pair of distinct subsets will be linearly separable under the flow map.
    
    \item \textbf{Alternative activation functions}. \cref{thm2} applies broadly to any Lipschitz-continuous activation function $\sigma:\R\to\R$, as long as there exist $ -\infty\leq a<b\leq \infty $ such that $\sigma(z) = 0$ for all $ z \in (a,b) $ and $\sigma\neq 0$. The choice of $\sigma$ would only affect the control time $T>0$. Following the terminology of reference \cite{Qianxiao2022}, such a $\sigma$ would be referred to as a well function, and our proofs could be extended using an argument based on affine-invariance of system \cref{eq:node}.
    
    Generalizing \cref{thm1} to other activations presents additional challenges since they must represent any flow of the truncated ReLU (as in \cref{lem:trunintorelu}). A modification (and improvement) of \cref{thm1} is \cref{thm:nodealt}, where we consider the triangular activation function.  In fact, any activation function $\sigma$ with $\operatorname{supp}(\sigma)\subset(a,b)$ for some $-\infty < a < b < +\infty$ allows for the same extension.

    \item \textbf{Alternative loss functions}. While we have used a simplified error function---the sum of the distances of all points to their respective target regions---cross-entropy loss $\ell_{\rm{CE}}$ is the standard choice for classification problems in practice. In binary classification, the cross entropy loss associated with a data point $(\bfx_n,y_n)\in\mathscr{D}$ is defined as \begin{equation*}
    \ell_{\rm{CE}}(\bfx_n,y_n)=-y_n\log(\hat y_n)-(1-y_n)\log(1-\hat y_n),
\end{equation*}
where $\hat y_n=\left(\operatorname{softmax}\circ P\circ\Phi_T(\bfx_n)\right)^{(1)}$ is the predicted probability that $y_n=1$. This prediction is obtained by applying a linear transformation $P:\R^d\to\R^2$ and then performing component-wise normalization using the softmax function. 

For any fixed $P$, our methods can be adapted to minimize this loss while maintaining the same complexity $L$. For example, if $d = 2$ and $P$ is the identity matrix, then we can achieve $\ell_{\rm{CE}}(\bfx_n, y_n) \to 0$ when $(-1)^{y_n} \, \Phi_T(\bfx_n) \cdot (-1, 1) \to +\infty$. To accomplish this, we can set $x^{(1)} = x^{(2)}$ as the decision boundary and increase the time $T$; however, this would not require additional switches. 
 \begin{figure}[t]  
\centering
\includegraphics[width=0.24\textwidth]{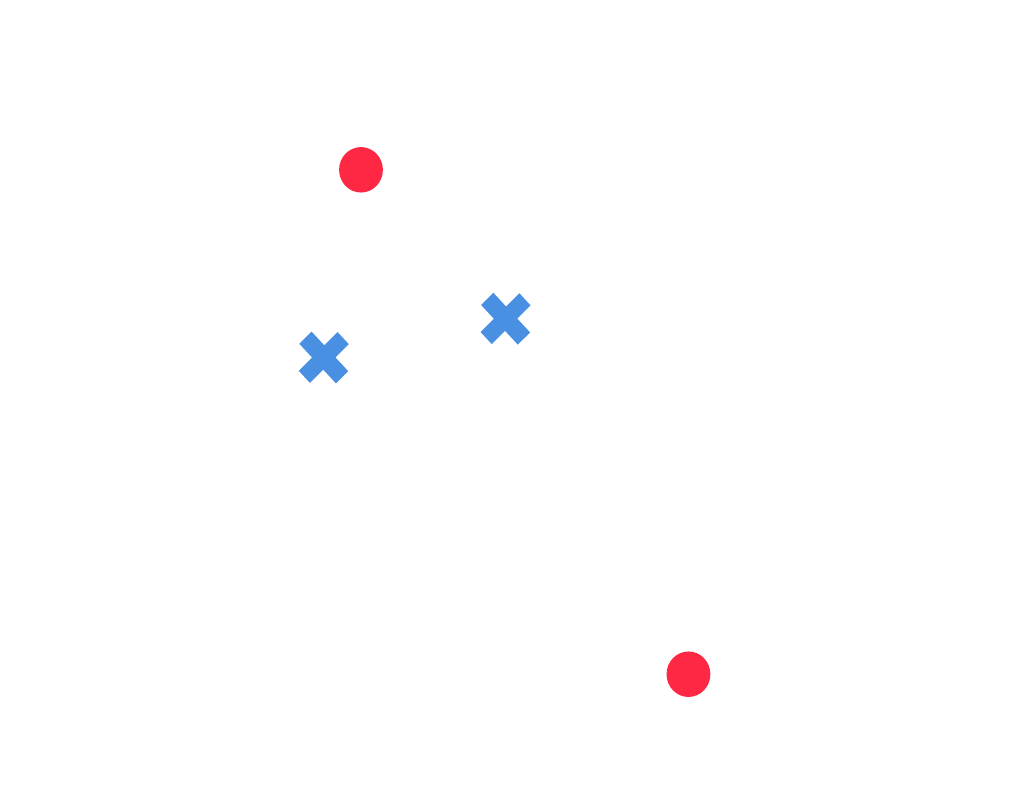}
\includegraphics[width=0.24\textwidth]{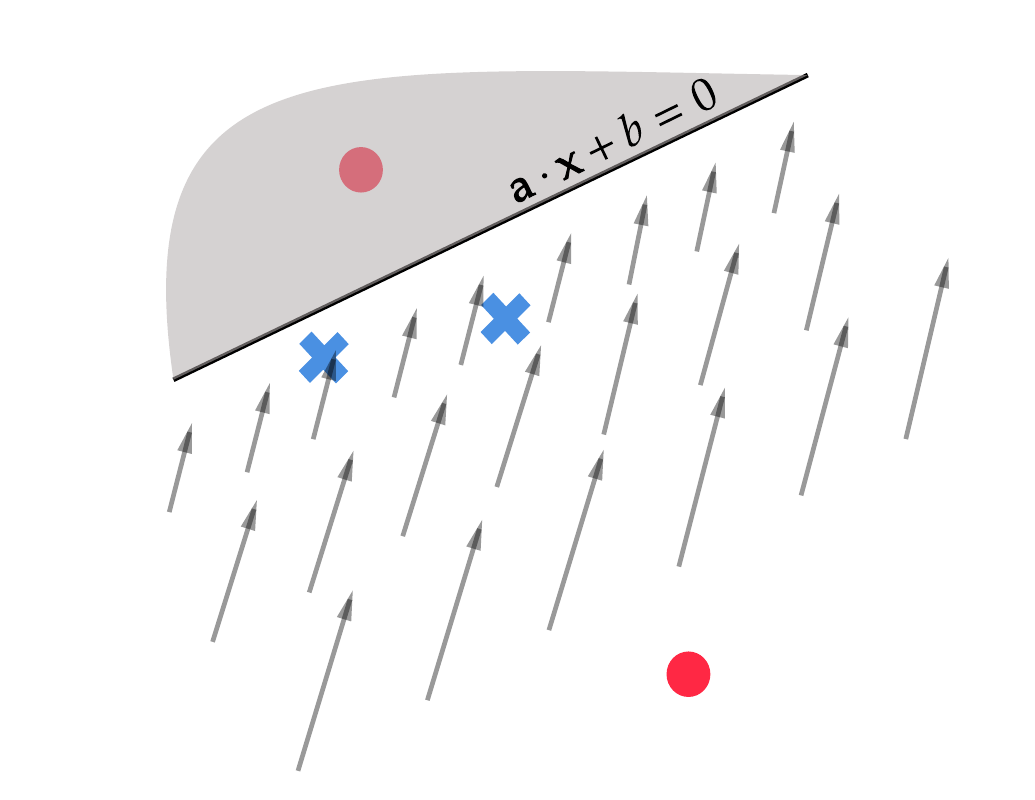}
\includegraphics[width=0.24\textwidth]{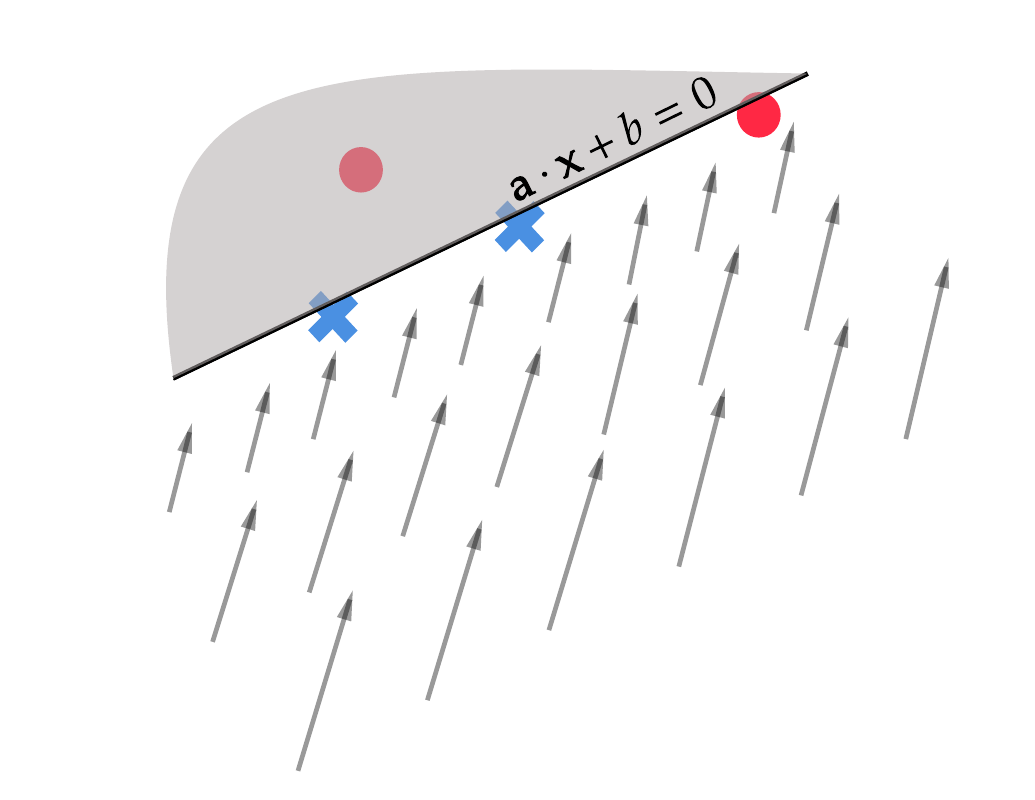}
\includegraphics[width=0.24\textwidth]{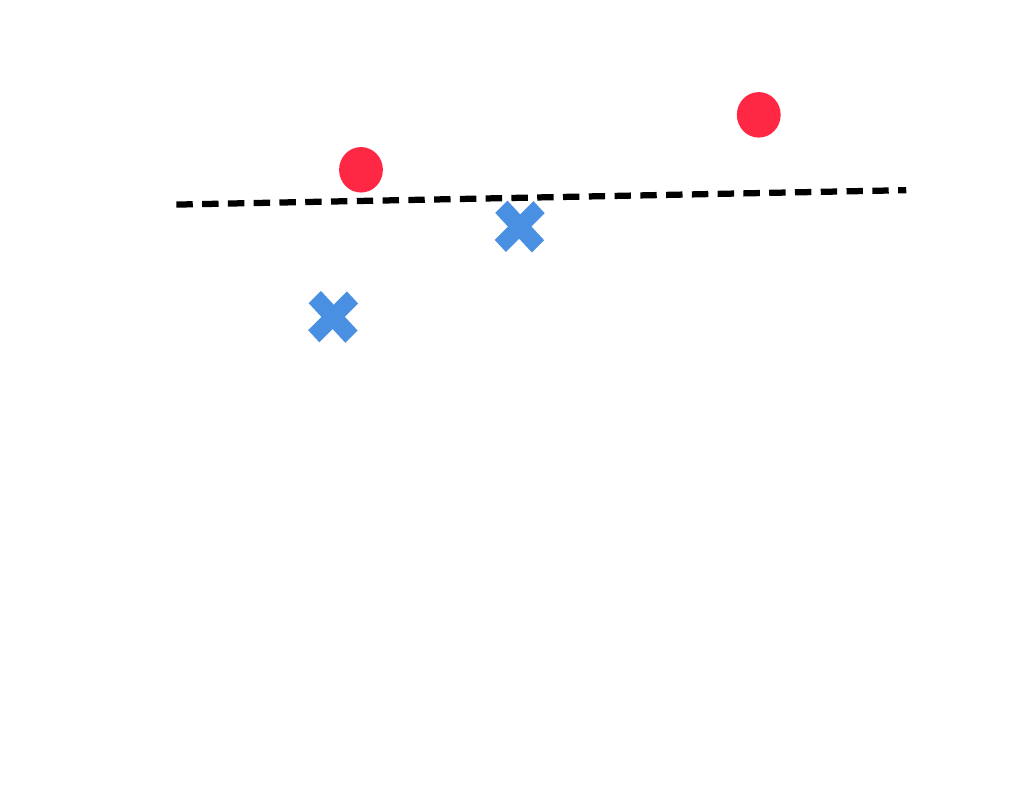}
\caption{In the first picture, at least two hyperplanes are required to separate both classes (``static'' classification). In the other three pictures, the same points evolve according to $\dot{\bfx} = \bfw (\bfa\cdot\bfx+b)_+$ (with constant $(\bfw,\bfa,b)$) and can eventually be separated using only one hyperplane (``dynamic'' classification).}
\label{fig:dynclas}
\end{figure}

\item \textbf{Point separability}. Since our control methods fundamentally rely on identifying hyperplanes that separate two classes (mostly via combinatorics), an important question arises: does dynamic classification using neural ODEs essentially reduce to static separability using hyperplanes? The short answer is \emph{no}; they are distinct frameworks.

On one hand, a collection of separating hyperplanes does not necessarily translate into effective neural ODE dynamics. In other words, having fixed $\bfa(t)$ and $b(t)$, there may be no control $\bfw(t)$ that yields the desired classification via \eqref{eq:node}. Our constructions require those hyperplanes to follow a specific structure (e.g., pairwise parallel in \cref{thm1}), making it possible to determine a suitable $\bfw(t)$. However, for a random family of separating hyperplanes, finding such $\bfw(t)$ can be highly complex. 

It might then seem that dynamic classification with neural ODEs is more restrictive than linear classification methods. Yet this is also false: neural ODEs can classify without explicitly relying on separating hyperplanes, that is, there exist controls $\bfa(t)$ and $b(t)$ that lead to successful classification even if the hyperplanes given by $\bfa(t)\cdot\bfx + b(t)=0$ do \emph{not} fully separate both classes. In \cref{fig:dynclas} we illustrate a scenario where neural ODEs outperform linear methods.

Our results address these questions by proving that, in the worst-case scenario (i.e., for any arbitrary dataset), neural ODEs can perform at least as well as linear classifiers. Furthermore, our numerical experiments seem to confirm that, under such worst-case conditions, the performance of neural ODEs trained via gradient descent matches our theoretical findings.

\item \textbf{Discrete neural networks}. In the discrete setting, the capacity to control $N$ points in $\R^d$ with a given level of complexity has been more broadly studied and often referred to as finite-sample expressivity. Early studies, focused on sigmoid activation functions, showed that networks with one hidden layer could perform this task using $O(N)$ neurons \cite{Huang91}.  With the rise of deep models and the ReLU activation, researchers explored how adding depth could reduce the required complexity \cite{powerofdepth, Huang2003}, and existent results were extended to ReLU networks, see \cite{understandingdl}. For deep feed-forward networks, \cite{Yun2018SmallRN} solved the problem using $O(\sqrt{N})$ neurons. 

In the past decade, the focus has shifted towards new architectures like convolutional networks \cite{nguyen_optimization_2018}, or ResNets \cite{hardt2017identity}, which achieve the goal using $O(N\log N)$ hidden nodes, assuming a minimum distance between points. This bound was later improved to $O(N/d)$ in \cite{Yun2018SmallRN}, under the assumption of general position. Notably, for conditional networks, the complexity can be reduced to  $O(\log{N})$ neurons  \cite{koyuncu2023memorization}. 

\end{itemize}

\section*{Future work}
Several important directions still need to be explored. Below, we present some of them.

\begin{enumerate}
    \item \textbf{Improvement or sharpness of the bound}. A first and natural question is whether the new bound we have introduced for the minimal number of switches $L$ required to classify any dataset in general position is sharp or can be further improved. 
    
    On one hand, it turns out that equality in \cref{eq:maxboundgenpos}  can be achieved in some simple cases, for example, when $d=N=2$, as illustrated in \cref{fig:gp1gp2}. However, our algorithm does not fully exploit the capacity of the parameter space because the piecewise constant controls $\bfa(t)$ and $\bfw(t)$ repeat certain values throughout the algorithm. Moreover, the algorithm is designed to iteratively classify clusters of $d$ points but does not consider the possibility of classifying larger clusters, unlike \cref{alg:thm1}.
    
    If our bound for $L$ is indeed sharp, a related problem  is to identify the geometry of those configurations where this bound is attained, similar to \cref{thm3}.

\item \textbf{Decision boundary}. In this study, we have restricted the decision boundary to be any hyperplane defined by the equation $x^{(i)}=1$ for $i\in\{1,\dots,d\}$. Our control methods, detailed in \cref{alg:thm1,alg:thm2}, can easily be adapted to every hyperplane in $\R^d$.  However, the classification problem may also be addressed using more complex hypersurfaces than hyperplanes or by treating the decision boundary itself as an optimizable parameter. This leads to the open question of determining the optimal decision boundary that minimizes complexity while maximizing training accuracy.


\item \textbf{Generalization}. The ultimate goal of supervised learning models is to make accurate predictions on new, unseen data that were not part of the training dataset. While our methods aim to determine the minimum complexity required for a model to achieve zero training error, this often leads to overfitting and poor performance on test data. 

Proper generalization requires to model to capture the geometric patterns in the configuration of data points, rather than perfectly classifying every individual point regardless of their distribution. We must mention, however, that recent studies 
suggest overfitting does not necessarily imply poor generalization \cite{2020bartlettoverfitting}. Understanding the relationship between both concepts remains an active area of research \cite{22yangjml}.

In practice, generalization is typically achieved by adding a regularization term to the loss function. In the context of simultaneous control, some outlier data points could be misclassified to reduce the control cost. Developing constructive control methods that intentionally sacrifice some accuracy on the training data to improve the model's ability to generalize remains an open question.

    \item \textbf{Combination with self-attention}. Transformers are deep architectures that alternate between attention and feed-forward layers, adding residual connections. A highly simplified continuous-time version is the interacting particle system given by:
    \begin{equation*}
        \dot \bfx_i(t) = \bfw(t)\sigma(\bfa(t)\cdot\bfx_i(t)+b(t))+\sum_{j=1}^N\frac{e^{\bfx_i(t)\cdot B(t)\bfx_j(t)}}{\sum_{\ell=1}^{N}e^{\bfx_i(t)\cdot B(t)\bfx_{\ell}(t)}}V(t)\bfx_j(t).
    \end{equation*}
    Here, the vector $\bfx_i(t)$ represents the position of the $i$-th particle (or data point) at time $t$, whereas $V(t),B(t)\in\R^{d\times d}$ are additional controls. The attention term induces a clustering effect among the data points, causing them to concentrate towards certain limiting configurations, as studied in \cite{geshkovski_emergence_2023}. An open question is to quantify the reduction in complexity achieved by controlling this mechanism to cluster points by classes as a preprocessing step for classification. A constructive control of continuous-time transformers has been carried out in \cite{geshkovski_measure_measure_2024} for the mean-field formulation of this system restricted to the sphere.
    
\end{enumerate}

\section*{Acknowledgments}
This work has been supported by the Madrid Government (Comunidad de Madrid, Spain) under the multiannual Agreement with UAM in the line for the Excellence of the University Research Staff in the context of the V PRICIT (Regional Programme of Research and Technological Innovation) and by PID2023-146872OB-I00 of MICIU (Spain). AÁ has been funded by FPU21/05673 from the Spanish Ministry of Universities.
RO  has been funded by  Severo
Ochoa Programme for Centres of Excellence in R\&D (CEX2023-001347-S) of MICIU (Spain).
EZ has been funded by the Alexander von Humboldt-Professorship program, the ModConFlex Marie Curie Action, HORIZON-MSCA-2021-DN-01, the COST Action MAT-DYN-NET, the Transregio 154 Project “Mathematical Modelling, Simulation and Optimization Using the Example of Gas Networks” of the DFG and TED2021-131390B-I00 of MICINN (Spain).

\bibliographystyle{unsrt}
\bibliography{biblio}

\begin{thebibliography}{10}

\bibitem{dlnature}
Yann LeCun, Y.~Bengio, and Geoffrey Hinton.
\newblock Deep learning.
\newblock {\em Nature}, 521:436--44, 05 2015.

\bibitem{imagelecun}
Y.~Lecun, L.~Bottou, Y.~Bengio, and P.~Haffner.
\newblock Gradient-based learning applied to document recognition.
\newblock {\em Proceedings of the IEEE}, 86(11):2278--2324, 1998.

\bibitem{languagebengio}
Yoshua Bengio, R\'{e}jean Ducharme, Pascal Vincent, and Christian Janvin.
\newblock A neural probabilistic language model.
\newblock {\em J. Mach. Learn. Res.}, 3:1137–1155, 2003.

\bibitem{fisher}
R.~A. Fisher.
\newblock The use of multiple measurements in taxonomic problems.
\newblock {\em Annals of Eugenics}, 7(7):179--188, 1936.

\bibitem{svm}
Corinna Cortes and Vladimir Vapnik.
\newblock Support-vector networks.
\newblock {\em Machine Learning}, 20(3):273--297, 1995.

\bibitem{randomforests}
Leo Breiman.
\newblock Random forests.
\newblock {\em Machine Learning}, 45(1):5--32, 2001.

\bibitem{rosenblatt1958perceptron}
Frank Rosenblatt.
\newblock The perceptron: a probabilistic model for information storage and
  organization in the brain.
\newblock {\em Psychological review}, 65(6):386, 1958.

\bibitem{weinan}
Weinan E.
\newblock A proposal on machine learning via dynamical systems.
\newblock {\em Communications in Mathematics and Statistics}, 5:1--11, 02 2017.

\bibitem{he2016residual}
Kaiming He, Xiangyu Zhang, Shaoqing Ren, and Jian Sun.
\newblock {Deep Residual Learning for Image Recognition}.
\newblock In {\em Proceedings of 2016 IEEE Conference on Computer Vision and
  Pattern Recognition}, pages 770--778, 2016.

\bibitem{ChenRBD18}
Ricky T.~Q. Chen, Yulia Rubanova, Jesse Bettencourt, and David Duvenaud.
\newblock Neural ordinary differential equations.
\newblock In {\em Proceedings of the 32nd International Conference on Neural
  Information Processing Systems}, page 6572–6583. Curran Associates Inc.,
  2018.

\bibitem{haber}
Eldad Haber and Lars Ruthotto.
\newblock Stable architectures for deep neural networks.
\newblock {\em Inverse Problems}, 34(1):014004, 2017.

\bibitem{massaroli2020dissecting}
Stefano Massaroli, Michael Poli, Jinkyoo Park, Atsushi Yamashita, and Hajime
  Asama.
\newblock Dissecting neural odes.
\newblock {\em Advances in Neural Information Processing Systems},
  33:3952--3963, 2020.

\bibitem{domzuazua}
Dom\`enec Ruiz-Balet and Enrique Zuazua.
\newblock Neural {ODE} control for classification, approximation, and
  transport.
\newblock {\em SIAM Rev.}, 65(3):735--773, 2023.

\bibitem{ALVAREZLOPEZ2024106640}
Antonio Álvarez{-}López, Arselane Hadj~Slimane, and Enrique Zuazua.
\newblock Interplay between depth and width for interpolation in neural odes.
\newblock {\em Neural Networks}, page 106640, 2024.

\bibitem{mumford1965geometric}
David Mumford, John Fogarty, and Frances Kirwan.
\newblock {\em Geometric Invariant Theory}.
\newblock Springer Berlin, 1994.

\bibitem{cover1965}
Thomas~M. Cover.
\newblock Geometrical and statistical properties of systems of linear
  inequalities with applications in pattern recognition.
\newblock {\em IEEE Transactions on Electronic Computers}, EC-14(3):326--334,
  1965.

\bibitem{cover_number_1967}
Thomas~M. Cover.
\newblock The {Number} of {Linearly} {Inducible} {Orderings} of {Points} in
  \textit{d} -{Space}.
\newblock {\em SIAM Journal on Applied Mathematics}, 15(2):434--439, March
  1967.

\bibitem{sontag_shattering_1997}
Eduardo~D. Sontag.
\newblock Shattering {All} {Sets} of ‘k’ {Points} in “{General}
  {Position}” {Requires} (k — 1)/2 {Parameters}.
\newblock {\em Neural Computation}, 9(2):337--348, February 1997.

\bibitem{TaGh}
Paulo Tabuada and Bahman Gharesifard.
\newblock Universal approximation power of deep residual neural networks
  through the lens of control.
\newblock {\em {IEEE} Trans. Autom. Control.}, 68(5):2715--2728, 2023.

\bibitem{CuLaTe}
Christa Cuchiero, Martin Larsson, and Josef Teichmann.
\newblock Deep neural networks, generic universal interpolation, and controlled
  {ODE}s.
\newblock {\em SIAM J. Math. Data Sci.}, 2(3):901--919, 2020.

\bibitem{EZOP}
Karthik Elamvazhuthi, Xuechen Zhang, Samet Oymak, and Fabio Pasqualetti.
\newblock Learning on manifolds: Universal approximations properties using
  geometric controllability conditions for neural odes.
\newblock In {\em Learning for Dynamics and Control Conference}, Proceedings of
  Machine Learning Research, pages 1--11, 2023.

\bibitem{Qianxiao2022}
Qianxiao Li, Ting Lin, and Zuowei Shen.
\newblock Deep learning via dynamical systems: An approximation perspective.
\newblock {\em Journal of the European Mathematical Society}, 25(5):1671--1709,
  2022.

\bibitem{cheng2023interpolation}
Jingpu Cheng, Qianxiao Li, Ting Lin, and Zuowei Shen.
\newblock Interpolation, approximation, and controllability of deep neural
  networks.
\newblock {\em SIAM Journal on Control and Optimization}, 63(1):625--649, 2025.

\bibitem{ishikawa_universal_2023}
Isao Ishikawa, Takeshi Teshima, Koichi Tojo, Kenta Oono, Masahiro Ikeda, and
  Masashi Sugiyama.
\newblock Universal {Approximation} {Property} of {Invertible} {Neural}
  {Networks}.
\newblock {\em Journal of Machine Learning Research}, 24(287):1--68, 2023.

\bibitem{uat1neuronresnet}
Hongzhou Lin and Stefanie Jegelka.
\newblock Resnet with one-neuron hidden layers is a universal approximator.
\newblock In {\em Proceedings of the 32nd International Conference on Neural
  Information Processing Systems}, page 6172–6181, 2018.

\bibitem{e_mean-field_2018}
Weinan E, Jiequn Han, and Qianxiao Li.
\newblock A mean-field optimal control formulation of deep learning.
\newblock {\em Research in the Mathematical Sciences}, 6(1):10, December 2018.

\bibitem{qianxiaomaxprinciple}
Qianxiao Li, Long Chen, Cheng Tai, and Weinan E.
\newblock Maximum principle based algorithms for deep learning.
\newblock {\em Journal of Machine Learning Research}, 18(165):1--29, 2018.

\bibitem{bonnet_measure_2023}
Benoît Bonnet, Cristina Cipriani, Massimo Fornasier, and Hui Huang.
\newblock A measure theoretical approach to the mean-field maximum principle
  for training {NeurODEs}.
\newblock {\em Nonlinear Analysis}, 227:113161, February 2023.

\bibitem{24isobejml}
Noboru Isobe and Mizuho Okumura.
\newblock Variational formulations of ode-net as a mean-field optimal control
  problem and existence results.
\newblock {\em Journal of Machine Learning}, 3(4):413--444, 2024.

\bibitem{Geshkovski_Zuazua_2022}
Borjan Geshkovski and Enrique Zuazua.
\newblock Turnpike in optimal control of pdes, resnets, and beyond.
\newblock {\em Acta Numerica}, 31:135–263, 2022.

\bibitem{EGe}
Carlos Esteve-Yag\"{u}e and Borjan Geshkovski.
\newblock Sparsity in long-time control of neural {ODE}s.
\newblock {\em Systems Control Lett.}, 172:Paper No. 105452, 14, 2023.

\bibitem{momentum}
Dom{\`e}nec Ruiz-Balet, Elisa Affili, and Enrique Zuazua.
\newblock Interpolation and approximation via momentum resnets and neural odes.
\newblock {\em Systems \& Control Letters}, 162:105182, 2022.

\bibitem{domzuazuaNormFlows}
Domènec Ruiz-Balet and Enrique Zuazua.
\newblock Control of neural transport for normalising flows.
\newblock {\em Journal de Mathématiques Pures et Appliquées}, 181:58--90,
  2024.

\bibitem{duchksep}
Wlodzislaw Duch.
\newblock K-separability.
\newblock In {\em Artificial Neural Networks - ICANN 2006}, volume 4131 of {\em
  Lecture Notes in Computer Science}, pages 188--197, 2006.

\bibitem{BoU1995}
Ralph~P. Boland and Jorge Urrutia.
\newblock Separating collections of points in {E}uclidean spaces.
\newblock {\em Inform. Process. Lett.}, 53(4):177--183, 1995.

\bibitem{freimer_complexity_1991}
Robert Freimer, Joseph S.~B. Mitchell, and Christine Piatko.
\newblock On the {Complexity} of {Shattering} {Using} {Arrangements}.
\newblock Technical {Report}, Cornell University, USA, March 1991.

\bibitem{AHMSS}
Esther~M. Arkin, Ferran Hurtado, Joseph S.~B. Mitchell, Carlos Seara, and
  Steven~S. Skiena.
\newblock Some lower bounds on geometric separability problems.
\newblock {\em Internat. J. Comput. Geom. Appl.}, 16(1):1--26, 2006.

\bibitem{HOULE1993139}
Michael~F. Houle.
\newblock Algorithms for weak and wide separation of sets.
\newblock {\em Discrete Applied Mathematics}, 45(2):139--159, 1993.

\bibitem{vc}
V.~N. Vapnik and A.~Ya. Chervonenkis.
\newblock On the uniform convergence of relative frequencies of events to their
  probabilities.
\newblock {\em Theory of Probability and its Applications}, 16(2):264--280,
  1971.

\bibitem{pytorch}
Adam Paszke, Sam Gross, Soumith Chintala, Gregory Chanan, Edward Yang, Zachary
  DeVito, Zeming Lin, Alban Desmaison, Luca Antiga, and Adam Lerer.
\newblock Automatic differentiation in pytorch.
\newblock In {\em NIPS 2017 Workshop on Autodiff}, 2017.

\bibitem{Kingma2014AdamAM}
Diederik~P. Kingma and Jimmy Ba.
\newblock Adam: A method for stochastic optimization.
\newblock {\em CoRR}, abs/1412.6980, 2014.

\bibitem{bless}
David Donoho.
\newblock High-dimensional data analysis: The curses and blessings of
  dimensionality.
\newblock {\em AMS Math Challenges Lecture}, pages 1--32, 01 2000.

\bibitem{Huang91}
S.-C. Huang and Y.-F. Huang.
\newblock Bounds on the number of hidden neurons in multilayer perceptrons.
\newblock {\em IEEE Transactions on Neural Networks}, 2(1):47--55, 1991.

\bibitem{powerofdepth}
Ronen {Eldan} and Ohad {Shamir}.
\newblock {The Power of Depth for Feedforward Neural Networks}.
\newblock {\em JMLR: Workshop and Conference Proceedings}, 49:1--34, 2015.

\bibitem{Huang2003}
Guang-Bin Huang.
\newblock Learning capability and storage capacity of two-hidden-layer
  feedforward networks.
\newblock {\em IEEE Transactions on Neural Networks}, 14(2):274--281, 2003.

\bibitem{understandingdl}
Chiyuan Zhang, Samy Bengio, Moritz Hardt, Benjamin Recht, and Oriol Vinyals.
\newblock Understanding deep learning requires rethinking generalization.
\newblock {\em Communications of the ACM}, 64, 11 2016.

\bibitem{Yun2018SmallRN}
Chulhee Yun, Suvrit Sra, and Ali Jadbabaie.
\newblock Small relu networks are powerful memorizers: A tight analysis of
  memorization capacity.
\newblock In {\em Proceedings of the 33rd International Conference on Neural
  Information Processing Systems}, page 15558–15569. Curran Associates Inc.,
  2019.

\bibitem{nguyen_optimization_2018}
Quynh Nguyen and Matthias Hein.
\newblock Optimization {Landscape} and {Expressivity} of {Deep} {CNNs}.
\newblock In {\em Proceedings of the 35th {International} {Conference} on
  {Machine} {Learning}}, pages 3730--3739. PMLR, July 2018.

\bibitem{hardt2017identity}
Moritz Hardt and Tengyu Ma.
\newblock Identity matters in deep learning.
\newblock In {\em International Conference on Learning Representations}, 2017.

\bibitem{koyuncu2023memorization}
Erdem Koyuncu.
\newblock Memorization capacity of neural networks with conditional
  computation.
\newblock In {\em The Eleventh International Conference on Learning
  Representations}, 2023.

\bibitem{2020bartlettoverfitting}
Peter~L. Bartlett, Philip~M. Long, Gábor Lugosi, and Alexander Tsigler.
\newblock Benign overfitting in linear regression.
\newblock {\em Proceedings of the National Academy of Sciences},
  117(48):30063--30070, 2020.

\bibitem{22yangjml}
Hongkang Yang.
\newblock A mathematical framework for learning probability distributions.
\newblock {\em Journal of Machine Learning}, 1(4):373--431, 2022.

\bibitem{geshkovski_emergence_2023}
Borjan Geshkovski, Cyril Letrouit, Yury Polyanskiy, and Philippe Rigollet.
\newblock The emergence of clusters in self-attention dynamics.
\newblock {\em Advances in Neural Information Processing Systems},
  36:57026--57037, December 2023.

\bibitem{geshkovski_measure_measure_2024}
Borjan Geshkovski, Philippe Rigollet, and Domènec Ruiz-Balet.
\newblock Measure-to-measure interpolation using {Transformers}, November 2024.
\newblock arXiv:2411.04551 [cs, math, stat].

\end{thebibliography}

\end{document}